\documentclass[11pt]{amsart}
\usepackage{amsmath,amssymb,amsthm}
\usepackage{version,tabularx,multicol,enumerate}
\usepackage{graphicx,float,psfrag}

\newcommand{\reff}[1]{(\ref{#1})}

\theoremstyle{plain}
\newtheorem{theo}{Theorem}[section]
\newtheorem{cor}[theo]{Corollary}
\newtheorem{prop}[theo]{Proposition}
\newtheorem{lem}[theo]{Lemma}
\newtheorem{defi}[theo]{Definition}
\theoremstyle{remark}
\newtheorem{rem}[theo]{Remark}

\newcommand{\AAke}{A_{S_k^\varepsilon}}

\newcommand{\ca}{{\mathcal A}}
\newcommand{\cb}{{\mathcal B}}

\newcommand{\cf}{{\mathcal F}}
\newcommand{\cg}{{\mathcal G}}
\newcommand{\ch}{{\mathcal H}}
\newcommand{\ci}{{\mathcal I}}
\newcommand{\cj}{{\mathcal J}}
\newcommand{\ck}{{\mathcal K}}
\newcommand{\cl}{{\mathcal L}}
\newcommand{\cn}{{\mathcal N}}
\newcommand{\cm}{{\mathcal M}}

\newcommand{\cs}{{\mathcal S}}

\newcommand{\ct}{{\mathcal T}}

\newcommand{\cw}{{\mathcal W}}

\newcommand{\D}{{\mathbb D}}
\newcommand{\E}{{\mathbb E}}

\newcommand{\M}{{\mathbb M}}
\newcommand{\N}{{\mathbb N}}
\renewcommand{\P}{{\mathbb P}}

\newcommand{\Q}{{\mathbb Q}}
\newcommand{\R}{{\mathbb R}}
\renewcommand{\S}{{\mathbb S}}

\newcommand{\rP}{{\rm P}}

\newcommand{\ind}{{\bf 1}}

\newcommand{\supp}{{\rm supp}\;}
\newcommand{\Supp}{{\rm Supp}\;}

\newcommand{\val}[1]{\mathop{\left| #1 \right|}\nolimits}
\newcommand{\inv}[1]{\mathop{\frac{1}{ #1}}\nolimits}
\newcommand{\expp}[1]{\mathop {\mathrm{e}^{ #1}}}

\newcommand{\lb}{[\![}
\newcommand{\rb}{]\!]}

\begin{document}

\includeversion{commentaries}

\title[Pruning a L\'evy CRT]{Pruning a L\'evy continuum random tree}

\date{\today}
\author{Romain Abraham} 

\address{
Romain Abraham,
MAPMO, CNRS UMR 6628,
F\'ed\'eration Denis Poisson FR 2964,
Universit\'e d'Orl\'eans,
B.P. 6759,
45067 Orl\'eans cedex 2
FRANCE.
}
  
\email{romain.abraham@univ-orleans.fr} 

\author{Jean-François Delmas}

\address{
Jean-Fran\c cois Delmas,
Universit\'e Paris-Est, CERMICS,  6-8
av. Blaise Pascal, 
  Champs-sur-Marne, 77455 Marne La Vall\'ee, FRANCE.}

\email{delmas@cermics.enpc.fr}

\author{Guillaume Voisin}

\address{
Guillaume Voisin,
MAPMO CNRS UMR 6628, 
F\'ed\'eration Denis Poisson FR 2964,
Universit\'e d'Orl\'eans,
B.P. 6759,
45067 Orl\'eans cedex 2
FRANCE.}

\email{guillaume.voisin@univ-orleans.fr}

\thanks{This work is partially supported by the 'Agence Nationale de
  la Recheche', ANR-08-BLAN-0190.}

\begin{abstract}
Given a general critical or sub-critical branching mechanism, we define
a pruning procedure of the associated L\'evy continuum random tree. This
pruning procedure is defined by adding some marks on the tree, using
L\'evy snake techniques. We then prove that the resulting sub-tree after
pruning is
still a L\'evy continuum random tree. This last result is proved using
the exploration process that codes the CRT, a special Markov
property and martingale problems for exploration processes. We finally
give the joint law under the excursion measure of the lengths of the
excursions of the initial exploration process and the pruned one.
\end{abstract}



\maketitle

\section{Introduction}

Continuous state branching processes (CSBP) were first introduced by Jirina
\cite{j:sbpcss} and it is known since Lamperti \cite{l:lsbp} that these
processes are the scaling limits of Galton-Watson processes. They
hence model the evolution of a large population on a long time interval. The law
of such a process is characterized by the so-called branching
mechanism function $\psi$. We will be interested mainly in critical or
sub-critical CSBP. In those cases, the branching mechanism $\psi$ is given
by 
\begin{equation}
\label{eq:psi}
\psi(\lambda)=\alpha\lambda+\beta\lambda^2+\int_{(0,+\infty)}
\pi(d\ell)\left(\expp{-\lambda\ell}-1+\lambda\ell\right), 
\quad \lambda\geq 0, 
\end{equation}
 with  $\alpha\ge 0$, $\beta\geq  0$ and  the L\'evy  measure $\pi$  is a
positive   $\sigma$-finite   measure    on   $(0,+\infty)$   such   that
$\int_{(0,+\infty)}  (\ell\wedge  \ell^2)\pi(d\ell)<\infty$. We shall
say that the branching mechanism $\psi$ has parameter $(\alpha,\beta,
\pi)$. Let us recall that $\alpha$ represents a drift term, $\beta$ is
a diffusion coefficient and $\pi$ describes the jumps of the CSBP.

As for discrete Galton-Watson processes, we can associate with a CSBP
a genealogical tree, see \cite{lglj:bplpep} or
\cite{dlg:rtlpsbp}. These trees can be considered as continuum
random trees (CRT) in the sense that the branching points along a
branch form a dense subset. We call the genealogical tree associated
with a branching mechanism $\psi$ the $\psi$-L\'evy CRT (the term
``L\'evy'' will be explained later). The prototype of such a tree is the
Brownian CRT introduced by Aldous \cite{a:crt1}. 

In a discrete setting, it is easy to consider and study a percolation
on the tree (for instance, see \cite{ap:tvmcdgwp} for percolation on the
branches of a Galton-Watson tree, or \cite{adh:pgwttxmp} for
percolation on the nodes of a Galton-Watson tree). 
The goal of this
paper is to introduce a general pruning procedure of a genealogical
tree associated with a branching mechanism $\psi$ of the form
\reff{eq:psi}, which is the continuous analogue of the previous
percolation (although no link is actually made between both). We first
add some marks on the skeleton of the tree according to a Poisson
measure with intensity $\alpha_1 \lambda$ where $\lambda$ is the
length measure on the tree (see the definition of that measure
further) and $\alpha_1$ is a non-negative parameter. We next add some
marks on the nodes of infinite index of the tree: with such a node $s$
is associated a ``weight'' say $\Delta_s$ (see later for a formal
definition), each infinite node is then marked with probability
$p(\Delta_s)$ where $p$ is a non-negative measurable function satisfying the
integrability condition
\begin{equation}\label{eq:cond-p}
\int_{(0,+\infty)}\ell\, p(\ell)\, \pi(d\ell)<+\infty.
\end{equation}
We then prune the tree according to these
marks and consider the law of the pruned subtree containing the root.
The main result of the paper is the following theorem:

\begin{theo}\label{thm:pruned_CRT}
Let $\psi$ be a (sub)-critical branching mechanism of the form
\reff{eq:psi}. We define
\begin{align}
\label{eq:def-pi0}d\pi_0(x) & :=\bigl(1-p(x)\bigr)d\pi(x)\\
\label{eq:def-alpha0}\alpha_0 & :=\alpha+\alpha_1+\int_{(,+\infty)}\ell p(\ell)\pi(d\ell)
\end{align}
and set
\begin{equation}\label{eq:def-psi0}
\psi_0(\lambda)=\alpha_0\lambda+\beta\lambda+\int_{(0,+\infty)}\pi_0(d\ell)\left(\expp{-\lambda\ell}-1+\lambda\ell\right)
\end{equation}
which is again a branching mechanism of a critical or subcritical
CSBP.

Then, the pruned subtree is a L\'evy-CRT with branching mechanism $\psi_0$.
\end{theo}

In order to make the previous statement more rigorous,  we must
first describe more precisely the geometric structure of a continuum
random tree and define the so-called exploration process that codes
the CRT in the next subsection.  
In a second subsection, we describe the pruning procedure  and  state
rigorously the  main  results of the paper. Eventually, we give some biological
motivations for studying the pruning procedure and other applications
of this work.

\subsection{The L\'evy CRT and its coding by the exploration process} 

We first give the definition of  a real tree, see e.g. \cite{e:prt} or
\cite{lg:rrt}.

\begin{defi}
A metric space $(\ct,d)$ is a real tree if the following two
properties hold for every $v_1,v_2\in\ct$.
\begin{itemize}
\item[(i)] There is a unique isometric map $f_{v_1,v_2}$
  from $[0,d(v_1,v_2)]$ into $\ct$ such that
$$f_{v_1,v_2}(0)=v_1\qquad\mbox{and}\qquad
  f_{v_1,v_2}(d(v_1,v_2))=v_2.$$
\item[(ii)] If $q$ is a continuous injective map from $[0,1]$ into
  $\ct$ such that $q(0)=v_1$ and $q(1)=v_2$, then we have
$$q([0,1])=f_{v_1,v_2}([0,d(v_1,v_2)]).$$
\end{itemize}
A rooted real tree is a real tree $(\ct,d)$ with a distinguished
vertex $v_\emptyset$ called the root.
\end{defi}

Let $(\ct,d)$ be a rooted real tree. The range of the mapping
$f_{v_1,v_2}$ is denoted by $\lb v_1,v_2,\rb$ (this is the line
  between $v_1$ and $v_2$ in the tree). In particular, for every
  vertex $v\in\ct$, $\lb v_\emptyset,v\rb$ is the path going from the
  root to $v$ which we call the ancestral line of vertex $v$. More
  generally, we say that a vertex $v$ is an ancestor of a vertex $v'$
  if $v\in\lb v_\emptyset,v'\rb$. If $v,v'\in \ct$, there is a unique
  $a\in \ct$ such that $\lb
  v_\emptyset,v\rb\cap\lb v_\emptyset,v'\rb=\lb v_\emptyset, a\rb$. We
  call $a$ the most recent common ancestor of $v$ and $v'$. By
  definition, the degree of a vertex $v\in\ct$ is the number of
  connected components of $\ct\setminus\{v\}$. A vertex $v$ is called
  a leaf if it has degree 1. Finally, we set $\lambda$ the
  one-dimensional Hausdorff measure on $\ct$.

The coding of a compact real tree by a continuous function is
now well known and is a key tool for defining random real trees. We
consider a continuous function $g\,:\, [0,+\infty)\longrightarrow
  [0,+\infty)$ with compact support and such that $g(0)=0$. We also
    assume that $g$ is not identically 0. For every $0\le s\le t$, we set
$$m_g(s,t)=\inf_{u\in[s,t]}g(u),$$
and
$$d_g(s,t)=g(s)+g(t)-2m_g(s,t).$$
We  then introduce the  equivalence relation  $s\sim t$  if and  only if
$d_g(s,t)=0$. Let  $\ct_g$ be the quotient  space $[0,+\infty)/\sim$. It
is easy  to check  that $d_g$ induces  a distance on  $\ct_g$. Moreover,
$(\ct_g,d_g)$  is a  compact  real tree  (see \cite{dlg:pfalt},  Theorem
2.1).  The function  $g$ is  the so-called  height process  of  the tree
$\ct_g$. This  construction can be  extended to more  general measurable
functions.

  In  order to define  a random  tree, instead  of taking  a tree-valued
  random variable,  it suffices  to take a  stochastic process  for $g$.
  For  instance,  when  $g$  is  a normalized  Brownian  excursion,  the
  associated real tree is Aldous' CRT \cite{a:crt3}.

The construction of a height process that codes a tree associated with
a general branching mechanism is due to Le Gall and Le Jan \cite{lglj:bplpep}.
Let $\psi$ be a branching mechanism given by \reff{eq:psi} and let $X$
be a L\'evy process with Laplace exponent $\psi$: $\E[\expp{-\lambda
    X_t}]=\expp{t\psi(\lambda)}$ for all $\lambda\ge 0$.  Following
\cite{lglj:bplpep}, we also  assume  that  $X$ is  of
infinite   variation    a.s.    which   implies    that   $\beta>0$   or
$\int_{(0,1)}\ell\pi(d\ell)=\infty$.    Notice that  these conditions   are
satisfied in the stable case: $\psi(\lambda)=\lambda^c$, $c\in
(1,2]$ (the quadratic case $\psi(\lambda)=\lambda^2$ corresponds to
  the Brownian case).

We then set
\begin{equation}\label{eq:defHbis}
H_t=\liminf_{\varepsilon\to
  0}\frac{1}{\varepsilon}\int_0^t\ind_{\{X_s<I_t^s+\varepsilon\}}ds
\end{equation}
where  for $0\le s\le t$,
$I_t^s=\inf_{s\le r \le t}X_r$.
If the additional assumption
\begin{equation}\label{eq:condHcont}
\int_1^{+\infty}\frac{du}{\psi(u)}<\infty
\end{equation}
holds, then the  process $H$ admits a continuous  version. In this case,
we  can consider  the  real tree  associated  with an  excursion of  the
process $H$ and we say that  this real tree is the L\'evy CRT associated
with $\psi$. If we set $L_t^a(H)$ the local time time of the process $H$
at level  $a$ and time  $t$ and $T_x=\inf\{t\ge 0,\  L_t^0(H)=x\}$, then
the process  $(L^a_{T_x}(H),a\ge 0)$  is a CSBP  starting from  $x$ with
branching mechanism $\psi$  and the tree with height  process $H$ can be
viewed as  the genealogical tree  of this CSBP.  Let us remark  that the
latter   property  also  holds   for  a   discontinuous  $H$   (i.e.  if
\reff{eq:condHcont} doesn't  hold) and we  still say that  $H$ describes
the genealogy of the CSBP associated with $\psi$.

In general, the process $H$ is not a Markov process. So, we introduce
the so-called exploration process $\rho=(\rho_t,t\ge 0)$ which is a
measure-valued process defined by
\begin{equation}\label{eq:def_rho}
\rho_t(dr)=\beta\ind_{[0,H_t]}(r)\; dr + \sum_{\stackrel{0<s\le t}
  {X_{s-}<I_t^s}}(I_t^s-X_{s-})\delta_{H_s}(dr). 
\end{equation}
 The  height process  can  easily  be  recovered from  the
exploration  process  as  $H_t=H(\rho_t)$,  where $H(\mu)$  denotes  the
supremum of the closed support of the measure $\mu$ (with the convention
that $H(0)=0$). If we endow the set $\cm_f(\R_+)$ of finite measures
on $\R_+$ with the topology of weak convergence, then 
the exploration process $\rho$ is a c\`ad-l\`ag
strong Markov process in $\cm_f(\R_+)$ (see
\cite{dlg:rtlpsbp}, Proposition 1.2.3).

To understand the meaning of the exploration process, let us use the queuing
system representation of \cite{lglj:bplpep} when $\beta=0$.  We consider
a preemptive LIFO (Last In, First Out) queue with one server.  A jump of
$X$ at time $s$ corresponds to the arrival of a new customer requiring a
service  equal  to $\Delta_s:=X_s-X_{s-}$.   The  server interrupts  his
current  job and  starts immediately  the service  of this  new customer
(preemptive  LIFO procedure).  When  this new  service is  finished, the
server  will  resume the  previous  job.  When $\pi$  is
infinite, all services will suffer interruptions.  The customer (arrived
at time) $s$  will still be in the  system at time $t>s$ if  and only if
$\displaystyle  X_{s-}<\inf_{s\le r\le  t}X_r$  and, in  this case,  the
quantity $\rho_t(\{H_s\})$ represents the  remaining service required by the
customer $s$ at time  $t$. Observe that $\rho_t([0,H_t])$ corresponds to
the load of the server at time $t$ and is equal to $X_t-I_t$ where
$$I_t=\inf\{X_u, 0\le u\le t\}.$$

In view of the Markov property of $\rho$ 
and  the Poisson  representation of  Lemma \ref{lem:dlg-decomp},  we can
view  $\rho_t$  as  a  measure  placed  on the  ancestral  line  of  the
individual labeled  by $t$ which  gives the intensity of  the sub-trees
that are grafted ``on the right'' of this ancestral line. The continuous
part of the measure $\rho_t$  gives binary branching points (i.e. vertex
in the tree of degree 3) which are dense along that ancestral line since
the excursion measure $\N$ that appears in Lemma \ref{lem:dlg-decomp} is
an infinite  measure, whereas  the atomic part  of the  measure $\rho_t$
gives nodes of infinite degree for the same reason.

Consequently, the nodes of the tree coded by $H$ are of two types :
nodes of degree 3 and nodes of infinite degree. Moreover, we see that
each node of infinite degree corresponds to a jump of the L\'evy
process $X$ and so we associate to such a node  a ``weight'' given by
the  height of  the corresponding  jump of  $X$ (this  will  be formally
stated in Section  \ref{sec:LPS}).  From now-on, we will  only handle the
exploration process although we will often use vocabulary taken from the
real  tree  (coded by  this  exploration  process).  In particular,  the
theorems will  be stated  in terms of  the exploration process  and also
hold when $H$ is not continuous.

\subsection{The pruned exploration process}

We now consider the L\'evy CRT associated with a general critical or
sub-critical branching mechanism $\psi$ (or rather the exploration
process that codes that tree) and we add marks on the tree. There will
be two kinds of marks: some marks will be set only on nodes of
infinite degrees whereas the others will be 'uniformly distributed' on
the skeleton on the tree.

\subsubsection{Marks on the nodes}

Let $p:[0,+\infty)\longrightarrow [0,1]$ be a measurable function
  satisfying condition \reff{eq:cond-p}.
Recall that each node of infinite degree of the tree is associated
with a jump $\Delta _s$ of the process $X$. Conditionally on $X$, we mark such a node with
probability $p(\Delta_s)$, independently of the other nodes.

\subsubsection{Marks on the skeleton}

Let $\alpha_1$ be a non-negative constant. The marks associated with
these parameters will be distributed on the skeleton of the tree
according to a Poisson point measure with intensity
$\alpha_1\lambda(dr)$ (recall that $\lambda$ denotes the
one-dimensional Hausdorff measure on the tree).

\subsubsection{The marked exploration process}

As we don't use the real trees framework but only the exploration
processes that codes the L\'evy CRTs, we must describe all these marks
in term of exploration processes. Therefore, we define a
measure-valued process
$$\cs:=((\rho_t,m^{\text{nod}}_t,m^{\text{ske}}_t),t\ge   0)$$
called the marked exploration process where the process $\rho$ is the usual
exploration process whereas the processes $m^{\text{nod}}$ and
$m^{\text{ske}}$ keep track of the marks, respectively on the nodes
and on the skeleton of the tree.

The measure $m^{\text{nod}}_t$ is just the sum of the Dirac measure of
the marked nodes (up to some weights for technical reasons) which are
the ancestors of $t$.

To define the measure $m_t^{\text{ske}}$, we first consider a L\'evy
snake $(\rho_t,W_t)_{t\ge 0}$ with spatial motion $W$ a Poisson
process of parameter $\alpha_1$ (see \cite{dlg:rtlpsbp}, Chapter 4 for
the definition of a L\'evy snake). We then define the measure
$m_t^{\text{ske}}$ as the derivative of the function $W_t$. Let us
remark that in \cite{dlg:rtlpsbp}, the height process is supposed to be
continuous for the construction of L\'evy snakes. We explain in the
appendix how to remove this technical assumption.

\subsubsection{Main result}

We denote by $A_t$ the Lebesgue measure of the set of the individuals
prior to $t$ whose lineage does not contain any mark i.e.
$$A_t=\int_0^t\ind_{\{m_s^{\text{nod}}=0,m_s^{\text{ske}}=0\}}ds.$$
We consider its right-continuous inverse $C_t:=\inf\{r\ge 0,\ A_r>t\}$
and we define the pruned exploration process $\tilde \rho$ by
$$\forall t\ge 0,\quad \tilde \rho_t=\rho_{C_t}.$$

In other  words, we remove from the  CRT all the individuals  who have a
marked  ancestor, and the  exploration process  $\tilde \rho$  codes the
remaining  tree.

We can now restate Theorem \ref{thm:pruned_CRT} rigorously in terms
of exploration processes.
\begin{theo}\label{thm:law_pruned_bref}
The  pruned exploration process $\tilde
\rho$ is distributed as  the exploration process
associated with a L\'evy process with Laplace exponent
$\psi_0$.
\end{theo}

The proof relies on a martingale problem for $\tilde \rho$ and a special
Markov  property, Theorem \ref{th:SMP}.   Roughly speaking,  the special
Markov property  gives the  conditional distribution of  the individuals
with marked  ancestors with respect to  the tree of  individuals with no
marked ancestors.  This result is  of independent interest.   Notice the
proof of this  result in the general setting  is surprisingly much more
involved than the previous two particular cases: the quadratic case (see
Proposition 6  in \cite{as:psf}  or Proposition 7  in \cite{blglj:sbps})
and   the   case  without   quadratic   term   (see   Theorem  3.12   in
\cite{ad:falp}).

Finally, we  give the joint  law of the  length of the excursion  of the
exploration  process and  the  length  of the  excursion  of the  pruned
exploration process, see Proposition \ref{prop:s-sq}.

\subsection{Motivations and applications}

A first approach for this construction is to consider the CSBP
$Y^0$ associated with the pruned exploration process $\tilde \rho$ as
an initial Eve-population which undergoes some neutral mutations (the
marks on the genealogical tree) and the CSBP $Y$ denotes the total
population (the Eve-one and the mutants) associated with the
exploration process $\rho$. We see that, from our construction, we
have
$$Y_0^0=Y_0,\qquad \mbox{and}\qquad \forall t\ge 0,\ Y_t^0\le Y_t.$$
The condition
$$d\pi_0(x)=(1-p(x))d\pi(x)$$
means that, when the population $Y^0$ jumps, so does the population
$Y$. By these remarks, we can see that our pruning procedure is quite
general. Let us however remark that the coefficient diffusion $\beta$
is the same for $\psi$ and $\psi_0$ which might imply that more
general prunings exist (in particular, we would like to remove some
of the vertices of index 3).

As  we  consider  general
critical or sub-critical branching mechanism, this work extends previous
work from Abraham and Serlet \cite{as:psf} on Brownian CRT ($\pi=0$) and
Abraham  and   Delmas  \cite{ad:falp}  on  CRT   without  Brownian  part
($\beta=0$). See also Bertoin  \cite{b:saptpgwpnm} for an approach using
Galton-Watson trees  and $p=0$, or \cite{ad:cbpimcsbpi} for an
approach using CSBP with immigration. Let  us remark that this  paper goes
along the  same general  ideas as \cite{ad:falp}  (the theorems  and the
intermediate lemmas  are the same)  but the proofs  of each of  them are
more involved and use quite different techniques based on martingale
problem.

This work  has also others  applications.  Our method separates  in fact
the genealogical  tree associated with $Y$ into  several components. For
some values of the parameters of the pruning procedure, we can construct
via our pruning procedure, a fragmentation process as defined by Bertoin
\cite{b:rfcp}  but   which  is   not  self-similar,  see   for  instance
\cite{as:psf, ad:falp,  v:dmfglt}.  On the  other hand, we can  view our
method as a manner to increase the size of a tree, starting from the CRT
associated with $\psi_0$  to get the CRT associated  with $\psi$. We can
even construct a tree-valued process which makes the tree grow, starting
from  a   trivial  tree  containing   only  the  root  up   to  infinite
super-critical trees, see \cite{ad:ctvmp}.

\subsection{Organization of the paper}

We first recall in the  next Section the construction of the exploration
process, how  it codes a  CRT and its  main properties we shall  use. We
also define the marked exploration process  that is used for pruning the tree. In
Section  \ref{sec:pruned}, we define  rigorously the  pruned exploration
process      $\tilde\rho$     and     restate      precisely     Theorem
\ref{thm:law_pruned_bref}. The rest of the paper is devoted to the proof
of that theorem. In Section \ref{sec:Markov_special}, we state and prove
a special Markov  property of the marked exploration process,  that gives the law
of the exploration process ``above'' the marks, conditionally on $\tilde \rho$. We use
this special property in Section \ref{sec:law_pruned} to derive from the
martingale problem  satisfied by $\rho$,  introduced in \cite{ad:fpigep}
when $\beta=0$, a martingale problem for $\tilde\rho$ which allows us to
obtain the law of $\tilde\rho$. Finally, we compute in the last section,
under the  excursion measure,  the joint law  of 
the lengths of the excursions of $\rho$ and $\tilde\rho$. The Appendix is
devoted to some extension of the L\'evy snake when the height process is
not continuous. 


\section{The exploration process: notations and properties}
\label{sec:levysnake}

We  recall  here  the  construction   of  the  height  process  and  the
exploration  process that  codes a  L\'evy continuum  random  tree. These
objects  have been  introduced  in \cite{lglj:bplpep,lglj:bplplfss}  and
developed later  in \cite{dlg:rtlpsbp}. The results of  this section are
mainly extracted from \cite{dlg:rtlpsbp}, but for Section \ref{sec:LPS}.

We denote by $\R_+$ the set of non-negative real numbers. Let   $\cm(\R_+)$
(resp.  $\cm_f(\R_+)$)  be  the  set  of  $\sigma$-finite  (resp.   finite)
measures  on  $\R_+$, endowed  with  the  topology  of vague  (resp.  weak)
convergence. If $E$ is a Polish space, let $\cb(E)$  (resp. $\cb_+(E)$) be the  set of
real-valued  measurable (resp. and non-negative) functions  defined on $E$
endowed    with   its   Borel    $\sigma$-field. For any measure $\mu\in\cm(\R_+)$ and $f\in \cb_+(\R_+)$, we write
$$\langle \mu,f\rangle =\int f(x)\,\mu(dx).$$

\subsection{The underlying L\'evy process}\label{subsec:levy}
We consider a $\R$-valued  L\'evy process $X=(X_t,t\ge 0)$ starting from
0. We  assume that $X$  is the canonical  process on the  Skorohod space
$\D(\R_+,\R)$ of  c\`ad-l\`ag real-valued paths,  endowed with the
canonical filtration.  The law of the  process $X$
starting  from  0  will  be   denoted  by  $\P$  and  the  corresponding
expectation by $\E$. Most of the following facts on L\'evy processes can
be found in \cite{b:pl}.

In this paper, we assume that $X$ 
\begin{itemize}
\item  has no negative jumps, 
\item  has first moments, 
\item  is of infinite variation,
\item  does not drift to $+\infty$. 
\end{itemize}

\medskip
The law of $X$  is characterized by its Laplace transform:
\[
\forall \lambda\ge 0,\qquad \E\left[\expp{-\lambda
    X_t}\right]=\expp{t\psi(\lambda)}
\]
where, as  $X$ does not drift to
$+\infty$,  its   Laplace  exponent  $\psi$  can  then   be  written  as
\reff{eq:psi}, 
where the L\'evy measure $\pi$ is a Radon measure on $\R_+$ (positive
    jumps) that satisfies the integrability condition
$$\int_{(0,+\infty)}(\ell\wedge\ell^2)\pi(d\ell)<+\infty$$
($X$ has first moments), the drift coefficient $\alpha$ is non
    negative ($X$ does not drift to $+\infty$) and $\beta\ge 0$. As we
    ask for $X$ to be of infinite variation, we must additionally
    suppose that $\beta>0$ or
    $\int_{(0,1)}\ell\,\pi(d\ell)=+\infty$.

As  $X$ is   of infinite  variation, we  have, see
Corollary VII.5 in \cite{b:pl},
\begin{equation}
   \label{eq:psi/l}
\lim_{\lambda
  \rightarrow\infty } \frac{\lambda}{\psi(\lambda)}
=0.
\end{equation}

Let $I=(I_t,t\ge 0)$
be  the infimum  process of  $X$, $I_t=\inf_{0\le  s\le t}X_s$,  and let
$S=(S_t,t\ge 0)$  be the supremum process,  $S_t=\sup_{0\le s\le t}X_s$.
We will  also consider for every $0\le  s\le t$ the infimum  of $X$ over
$[s,t]$:
\[
I_t^s=\inf_{s\le r\le t}X_r.
\]

We denote by $\cj$ the set of jumping times of $X$:
\begin{equation}
   \label{eq:def-cj}
\cj=\{t\ge 0,\ X_t>X_{t-}\}
\end{equation}
and for $t\ge 0$ we set $\Delta _t:=X_t-X_{t-}$ the height of
the jump of $X$ at time $t$. Of course, $\Delta _t>0\iff t \in\cj$.

The point 0 is regular for the Markov 
process $X-I$, and $-I$ is the  local time of $X-I$ at 0 (see
\cite{b:pl}, chap. VII). Let $\N$ be the associated excursion measure of
the process $X-I$ away from 0, and let $\sigma=\inf\{t>0; X_t-I_t=0\}$
be the length
of the excursion of $X-I$ under $\N$. We will assume that under $\N$,
$X_0=I_0=0$. 

Since $X$ is of infinite variation, 0 is also regular for the Markov 
process $S-X$. The local time $L=(L_t, t\geq 0)$ of $S-X$ at 0 will be
normalized so that 
\[
\E[\expp{-\lambda S_{L^{-1}_t}}]= \expp{- t \psi(\lambda)/\lambda},
\]
where $L^{-1}_t=\inf\{ s\geq 0; L_s\geq t\}$ (see also \cite{b:pl}
Theorem VII.4 (ii)).

\subsection{The height process}

We now define the height process $H$ associated with the L\'evy
process $X$. Following \cite{dlg:rtlpsbp}, we give an alternative
definition of $H$ instead of those in the introduction, formula
\reff{eq:defHbis}.

For each $t\geq 0$, we consider the reversed process at time $t$,
$\hat X^{(t)}=(\hat X^{(t)}_s,0\le s\le t)$ by:
\[
\hat X^{(t)}_s=
X_t-X_{(t-s)-} \quad \mbox{if}\quad  0\le s<t,
\]
with the convention $X_{0-}=X_0$. The two processes $(\hat X^{(t)}_s,0\le s\le t)$
and $(X_s,0\le s\le t)$ have the same law. Let $\hat S^{(t)}$ be the
supremum process of $\hat X^{(t)}$ and $\hat L^ {(t)}$ be the
local time at $0$ of $\hat S^{(t)} - \hat X^{(t)}$ with the same
normalization as $L$.

\begin{defi}\label{def:height_process}(\cite{dlg:rtlpsbp}, Definition
  1.2.1)\\
  There  exists  a lower semi-continuous modification of the process $(\hat
  L^{(t)},t\ge 0)$. We denote by $(H_t,t\ge 0)$ this modification.
\end{defi}

This definition gives also a modification of the process defined by
\reff{eq:defHbis} (see \cite{dlg:rtlpsbp}, Lemma 1.1.3). In general,
$H$ takes its values in $[0,+\infty]$, but we have, a.s. for every
$t\ge 0$, $H_s<\infty$ for every $s<t$ such that $X_{s-}\le I_t^s$, and
$H_t<+\infty$ if $\Delta _t>0$ (see \cite{dlg:rtlpsbp}, Lemma
1.2.1).
The process $H$ does not admit a continuous version (or even
c\`ad-l\`ag) in general but it has continuous sample paths $\P$-a.s. iff
\reff{eq:condHcont} is satisfied, see \cite{dlg:rtlpsbp}, Theorem 1.4.3.

To end this section, let us remark that the height process is also
well-defined under the excursion process $\N$ and all the previous
results remain valid under $\N$.

\subsection{The exploration process}
\label{sec:PLRT}
The height process is not Markov in general. But it is a very simple
function of a measure-valued Markov process, the so-called exploration
process. 

The     exploration    process     $\rho=(\rho_t,t\ge    0)$     is    a
$\cm_f(\R_+)$-valued process defined  as follows: for  every $f\in
\cb_+(\R_+) $, $\langle \rho_t,f\rangle =\int_{[0,t]} d_sI_t^sf(H_s)$,
or equivalently 
\begin{equation}
\rho_t(dr)=\beta\ind_{[0,H_t]}(r)\; dr + \sum_{\stackrel{0<s\le t}
  {X_{s-}<I_t^s}}(I_t^s-X_{s-})\delta_{H_s}(dr). 
\end{equation}
In particular, the total mass of $\rho_t$ is
$\langle \rho_t,1\rangle =X_t-I_t$.

For $\mu\in \cm(\R_+)$, we set 
\begin{equation}
   \label{def:H}
H(\mu)=\sup\, \Supp \mu,
\end{equation}
where $ \Supp \mu$ is the closed support of $\mu$,  with the
convention $H(0)=0$. We have

\begin{prop}\label{prop:rho} (\cite{dlg:rtlpsbp}, Lemma 1.2.2 and
  Formula (1.12))\\
Almost surely, for every $t>0$,
\begin{itemize}
\item $H(\rho_t)=H_t$,
\item $\rho_t=0$ if and only if $H_t=0$,
\item if $\rho_t\neq 0$, then $\Supp \rho_t=[0,H_t]$. 
\item $\rho_t= \rho_{t^-} + \Delta_t \delta_{H_t}$, where $\Delta_t=0$
  if $t\not\in \cj$. 
\end{itemize}
\end{prop}

In the definition  of the exploration process, as $X$  starts from 0, we
have  $\rho_0=0$ a.s.  To  state the  Markov property  of $\rho$,  we must
first define  the  process  $\rho$  started  at any  initial  measure  $\mu\in
\cm_f(\R_+)$.

For $a\in  [0, \langle  \mu,1\rangle ] $,  we define the  erased measure
$k_a\mu$ by
\begin{equation}
   \label{eq:def-ka}
k_a\mu([0,r])=\mu([0,r])\wedge (\langle \mu,1\rangle -a), \quad \text{for $r\geq 0$}.
\end{equation}
If $a> \langle  \mu,1\rangle $, we set $k_a\mu=0$.   In other words, the
measure $k_a\mu$ is the measure $\mu$ erased by a mass $a$ backward from
$H(\mu)$.

For $\nu,\mu \in \cm_f(\R_+)$, and $\mu$ with compact support, we
  define the concatenation $[\mu,\nu]\in \cm_f(\R_+) $ of the
  two measures by: 
\[
\bigl\langle [\mu,\nu],f\bigr\rangle =\bigl\langle \mu,f\bigr\rangle +\bigl\langle \nu,f(H(\mu)+\cdot)\bigr\rangle ,
\quad f\in \cb_+(\R_+).
\]

Finally,  we  set for  every  $\mu\in  \cm_f(\R_+)$  and every  $t>0$,
\begin{equation}
   \label{eq:rhot-t-mu}
\rho_t^\mu=\bigl[k_{-I_t}\mu,\rho_t].
\end{equation}
  We say that $(\rho^\mu_t, t\geq
0)$  is  the  process  $\rho$  started at  $\rho_0^\mu=\mu$,  and  write
$\P_\mu$  for its  law. Unless  there is  an ambiguity,  we  shall write
$\rho_t$ for $\rho^\mu_t$.

\begin{prop}(\cite{dlg:rtlpsbp}, Proposition 1.2.3)\\
For any initial finite measure $\mu\in\cm_f(\R_+)$, the process
$(\rho_t^\mu,t\ge 0)$ is a c\`ad-l\`ag strong Markov process in 
$\cm_f(\R_+)$.
\end{prop}

\begin{rem}
  \label{rm:rho-L} 

{F}rom the construction of $\rho$, we get that a.s.
  $\rho_t=0$ if and only if $ -I_t\geq {\langle \rho_0,1\rangle }$ and
  $X_t-I_t=0$.  This implies that $0$ is also a regular point for
  $\rho$. 
  Notice that $\N$  is also the excursion measure  of the process $\rho$
  away  from $0$, and  that $\sigma$,  the length  of the  excursion, is
  $\N$-a.e.  equal to $\inf\{ t>0; \rho_t=0\}$.
\end{rem}

Exponential formula for the Poisson point process of jumps of the
inverse subordinator of $-I$ gives (see also the beginning of Section
3.2.2. \cite{dlg:rtlpsbp}) that for $\lambda>0$
\begin{equation}
   \label{eq:N_s}
\N\left[1 -\expp{-\lambda
  \sigma}\right] =\psi^{-1}(\lambda). 
\end{equation}

\subsection{The marked exploration process}
\label{sec:LPS}

As presented in the introduction, we add random marks on the L\'evy CRT
coded by $\rho$. There will be two kinds of marks: marks on the nodes
of infinite degree and marks on the skeleton.

\subsubsection{Marks  on   the  skeleton}
Let $\alpha_1\geq 0$. We want to construct a ``L\'evy Poisson snake''
(i.e. a L\'evy snake with spatial motion a Poisson process), whose  jumps
give the marks on the branches of the CRT. More precisely, we
set $\cw$ the space of killed c\`ad-l\`ag paths $w\,:\,
[0,\zeta)\rightarrow \R$ where $\zeta\in(0,+\infty)$ is called the
  lifetime of the path $w$. We equip $\cw$ with a distance $d$
  (defined in \cite{dlg:rtlpsbp} Chapter 4 and whose expression is not
  important for our purpose) such that $(\cw,d)$ is a Polish space.

  By Proposition 4.4.1 of  \cite{dlg:rtlpsbp} when $H$ is continuous, or
  the  results of  the  appendix in  the  general case,  there exists  a
  probability       measure      $\tilde       \P$       on      $\tilde
  \Omega=\D(\R_+,\cm_f(\R_+)\times  \cw)$   under  which  the  canonical
  process $(\rho_s,W_s)$ satisfies
\begin{enumerate}
\item The process $\rho$ is the exploration process starting at 0
  associated with a branching mechanism $\psi$,
\item For every $s\ge 0$, the path $W_s$ is distributed as a Poisson
  process with intensity $\alpha_1$ stopped at time $H_s:=H(\rho_s)$,
\item The process $(\rho,W)$ satisfies the so-called snake property:
  for every $s<s'$, conditionally given $\rho$, the paths $W_s(\cdot)$ and $W_{s'}(\cdot)$ coincide
  up to time $H_{s,s'}:=\inf\{H_u,s\le u\le s'\}$ and then are independent.
\end{enumerate}

So, for every $t\ge 0$, the path $W_t$ is a.s. c\`ad-l\`ag with jumps
equal to one. Its derivative $m_t^{\text{ske}}$ is an atomic measure
on $[0,H_t)$; it gives the marks (on the skeleton) on the ancestral line of the
  individual labeled $t$.

We shall denote by $\tilde \N$ the corresponding excursion measure out of
$(0,0)$. 

\subsubsection{Marks on the nodes}\label{sec:mark_nodes}
Let  $p$ be a  measurable function  defined on  $\R_+$ taking  values in
$[0,1]$ such that  
\begin{equation}\label{eq:condition-pi1}
\int_{(0,+\infty )} \ell p(\ell)\; \pi(d\ell)<\infty.
\end{equation}
We define the measures $\pi_1$ and $\pi_0$ by their density: 
$$d\pi_1(x)=p(x)d\pi(x)\qquad \mbox{and}\qquad
d\pi_0(x)=(1-p(x))d\pi(x).$$

Let $(\Omega',\ca',\rP')$  be a probability  space with no  atom. Recall
that $\cj$,  defined by \reff{eq:def-cj},  denotes the jumping  times of
the L\'evy process $X$ and that $\Delta _s$ represents the height of the
jump of $X$ at time $s\in \cj$.  As $\cj$ is countable, we can construct
on  the product  space  $\tilde \Omega\times\Omega'$  (with the  product
probability measure $\tilde \P\otimes  \rP'$) a family $(U_s, s\in \cj)$
of  random  variables  which  are, conditionally  on  $X$,  independent,
uniformly  distributed over $[0,1]$  and independent  of $(\Delta_s,s\in
\cj)$ and $(W_s, s\geq 0)$.  We set, for every $s\in \cj$:
\[
V_{s}=\ind_{\{U_s\le p(\Delta _{s})\}},
\]
so that,  conditionally on $X$, the
family $(V_s,s\in \cj)$ are independent Bernoulli random variables with
respective parameters $p(\Delta _s)$.

We  set $\cj^1=\{s\in\cj,\  V_s=1\}$ the  set  of the  marked jumps  and
$\cj^0=\cj\setminus\cj^1=\{s\in\cj,\ V_s=0\}$ the  set of the non-marked
jumps.  For $t\geq 0$, we consider the measure on $\R_+$,
\begin{equation}
   \label{eq:def-m-nod}
m^{\text{nod}}_t    (dr)   =   \sum_{\stackrel{0<s\le    t,\;s\in   \cj^1}
  {X_{s-}<I_t^s}}\left(I_t^s-X_{s-}\right)\delta_{H_s}(dr).
\end{equation}
The atoms of $m^{\text{nod}}_t$ give the marked nodes of the
exploration process at
time $t$.

The definition of the measure-valued process $m^{\text{nod}}$ also
holds under $\tilde \N\otimes \rP'$. For convenience, we shall write $\P$ for 
$\tilde \P\otimes \rP'$ and $\N$ for $\tilde \N\otimes \rP'$. 

\subsubsection{Decomposition of $X$}
At  this  stage,  we  can  introduce  a  decomposition  of  the  process
$X$.  Thanks to the  integrability condition  \reff{eq:condition-pi1} on
$p$, we can define the process $X^{(1)}$ by, for every $t\ge 0$,
$$X^{(1)}_t=\alpha_1 t+\sum_{0<s\le t;\ s\in\cj^1}\Delta_s.$$
The process $X^{(1)}$ is a subordinator with Laplace exponent $\phi_1$
given by:
\begin{equation}\label{eq:phi1}
\phi_1(\lambda)=\alpha_1\lambda+\int_{(0,+\infty)}
\pi_1(d\ell)\left(1-\expp{-\lambda\ell}\right), 
\end{equation}
with $\pi_1(dx)=p(x)\pi(dx)$.
We then set $X^{(0)}=X-X^{(1)}$ which is a L\'evy process with Laplace
exponent $\psi_0$, independent of the process $X^{(1)}$ by standard
properties of Poisson point processes.

We assume that $\phi_1\ne 0$ so that $\alpha_0$ defined by
\reff{eq:def-alpha0} is such that:
\begin{equation}\label{eq:alpha0>0}
\alpha_0>0.
\end{equation}
It is easy to check, using  $\int_{(0,\infty )}
\pi_1(d\ell)\ell<\infty $, that 
\begin{equation}
   \label{eq:lim-phi1}
\lim_{\lambda\rightarrow\infty } \frac{\phi_1(\lambda)}{\lambda}=\alpha_1.
\end{equation}

\subsubsection{The marked exploration process}

We consider the  process
$$\cs=((\rho_t,m_t^{\text{nod}},m_t^{\text{ske}}), t\geq 0)$$
on the  product probability space  $\tilde \Omega \times  \Omega'$ under
the probability $\P$ and call  it the marked exploration process. Let us
remark that,  as the process is  defined under the  probability $\P$, we
have $\rho_0=0$, $m_0^{\text{nod}}=0$ and $m_0^{\text{ske}}=0$ a.s.

Let us first define the state-space of the marked exploration process.
We consider the set
$\S$ of triplet $(\mu,\Pi_1,\Pi_2)$ where
\begin{itemize}
\item $\mu$ is a finite measure on $\R_+$,
\item $\Pi_1$ is a finite measure on $\R_+$ absolutely continuous with
  respect to $\mu$,
\item $\Pi_2$ is a $\sigma$-finite measure on $\R_+$ such that
\begin{itemize}
\item $\mathrm{Supp}(\Pi_2)\subset \mathrm{Supp}(\mu)$,
\item for every $x<H(\mu)$, $\Pi_2([0,x])<+\infty$,
\item if $\mu(\{H(\mu)\})>0$, $\Pi_2(\R_+)<+\infty$.
\end{itemize}
\end{itemize}

We endow $\S$ with the following distance: If $(\mu,\Pi_1,\Pi_2)\in\S$, we set
$$w(t)=\int \ind_{[0,t)}(\ell)\Pi_2(d\ell)$$
and
$$\tilde w(t)=w\left(H(k_{(\langle\mu,1\rangle -t)}\mu)\right)
\quad\text{for}\quad 
t\in [0, \langle\mu,1\rangle ).$$
We then define
$$d'((\mu,\Pi_1,\Pi_2),(\mu',\Pi_1',\Pi_2'))=d((\mu,\tilde
w),(\mu',\tilde w'))+D(\Pi_1,\Pi_1')$$
where $d$ is the distance defined by \reff{eq:dist_d} and $D$ is a
distance that defines the topology of weak convergence and such that
the metric space $(\cm_f(\R_+),D)$ is complete.

 To  get the Markov property of the marked
 exploration process, we  must define the process $\cs$  started at any initial
value of $ \S$.   For $(\mu,  \Pi^{\text{nod}}, \Pi^{\text{ske}})\in\S$,
we set     $\Pi=(\Pi^{\text{nod}},
\Pi^{\text{ske}} )$ and  $H^\mu_t=H(k_{-I_t}
\mu)$. We   define
$$(m^{\text{nod}})^{(\mu,  \Pi)}_t   =\left[\Pi^{\text{nod}}   \ind_{[0,
  H_t^\mu)}+\ind_{\{\mu(\{H_t^\mu\}) >0\}}\frac{k_{-I_t}\mu(\{H_t^\mu\})\Pi^{\text{nod}}(\{H_t^\mu\})}{\mu(\{H_t^\mu\})}\delta 
  _{H_t^\mu},
  m^{\text{nod}}_t\right]$$
and   
$$(m^{\text{ske}})^{(\mu,   \Pi)}_t  =[\Pi^{\text{ske}}   \ind_{[0,
  H_{t}^\mu)},    m^{\text{ske}}_t].$$
Notice    the    definition    of
$(m^{\text{ske}})^{(\mu, \Pi)}_t $ is  coherent with the construction of
the  L\'evy  snake,   with  $W_0$   being  the   cumulative   function  of
$\Pi^{\text{ske}}$ over $[0,H_0]$.

We shall write  $m^{\text{nod}}$ for $(m^{\text{nod}})^{(\mu, \Pi)}$ and
similarly      for     $m^{\text{ske}}$.      Finally,      we     write
$m=(m^{\text{nod}},m^{\text{ske}})$. By construction  and since $\rho$ is
an   homogeneous  Markov   process,  the   marked   exploration  process
$\cs=(\rho,m)$ is an homogeneous Markov process.

From now-on, we suppose that the marked exploration process is defined
on the canonical space $(\S,\cf')$ where $\cf'$ is the Borel
$\sigma$-field associated with the metric $d'$. We denote by
$\cs=(\rho,m^{\text{nod}},m^{\text{ske}})$ the canonical process and 
we denote by  $\P_{\mu,\Pi}$ the probability measure under which the
canonical process is distributed as  the  marked
exploration process
starting at time 0  from $(\mu,\Pi)$, and by $\P_{\mu,\Pi}^*$  the probability measure under which the
canonical process is distributed as  the
marked  exploration process  killed when  $\rho$ reaches  $0$. For  convenience
we shall  write $\P_\mu$ if $\Pi=0$ and $\P$ if $(\mu, \Pi)=0$ and
similarly for $\P^*$. Finally, we still denote by $\N$ the
distribution of $\cs$ when $\rho$ is distributed under the excursion
measure $\N$.

Let  $\cf=(\cf_t,  t\geq  0)$  be  the  canonical filtration.   Using  the  strong   Markov  property  of
$(X,X^{(1)})$  and  Proposition \ref{prop:annexe}  or  Theorem 4.1.2  in
\cite{dlg:rtlpsbp} if $H$ is continuous, we get the following result.
\begin{prop}
   The marked exploration process $\cs$ is a c\`ad-l\`ag $\S$-valued  strong Markov process. 
\end{prop}
Let us remark that the marked exploration process satisfies the
following snake property:
\begin{equation}\label{eq:snake-prop}
\P-a.s. \ (\mbox{or }\N-a.e.),\quad (\rho_t,m_t)(\cdot\cap
    [0,s])=(\rho_{t'},m_{t'})(\cdot\cap [0,s])\ \mbox{for every }0\le s<H_{t,t'}.
\end{equation}

\subsection{Poisson representation}\label{sec:poisson}

We  decompose the  path of  $\cs$ under  $\P^*_{\mu, \Pi}$  according to
excursions of  the total mass of  $\rho$ above its  past minimum, see Section
4.2.3 in  \cite{dlg:rtlpsbp}.  More precisely,   let
$(a_i,b_i), i\in \ck$ be the 
excursion  intervals of  $X-I$ above  $0$ under  $\P^*_{\mu,  \Pi}$. For
every $i\in \ck$, we define $h_i=H_{a_i}$  and
$\bar\cs^i=(\bar\rho^i, \bar m^i)$
by the formulas:  for $t\geq 0$ and $f\in \cb_+(\R_+)$, 
\begin{align}
\label{eq:rho-min}
\langle \bar\rho_t^i,f\rangle                                                             &
=\int_{(h_i,+\infty)}f(x-h_i)\rho_{(a_i+t)\wedge
b_i}(dx)\\                        
\label{eq:m-min}
\langle (\bar m_t^\text{a})^i,f\rangle                         &
=\int_{(h_i,+\infty)}f(x-h_i)m_{(a_i+t)\wedge
b_i}^\text{a}(dx),\quad \text{a}\in \{\text{nod},\text{ske}\},
\end{align}
with  $\bar m^i=((\bar m^\text{nod})^i, (\bar m^\text{ske})^i)$. We  set $\bar\sigma^i=\inf
\{s>0 ;  \langle \rho_s^i,1 \rangle =0  \}$.  It is easy  to adapt Lemma
4.2.4. of \cite{dlg:rtlpsbp} to get the following Lemma.
\begin{lem}
\label{lem:dlg-decomp}
   Let $(\mu, \Pi)\in\S$. The point
   measure $\displaystyle \sum_{i\in \ck} \delta_{(h_i,\bar\cs^i)}$ is under
    $\P^*_{\mu, \Pi}$ a Poisson point measure with intensity $\mu(dr)
    \N[d\cs]$. 
\end{lem}

\subsection{The dual process and representation formula}
\label{sec:dual}

We  shall need the  $\cm_f(\R_+)$-valued process  $\eta=(\eta_t,t\ge 0)$
 defined by
\[
\eta_t(dr)=\beta\ind_{[0,H_t]}(r)\; dr + \sum_{\stackrel{0<s\le
    t}{X_{s-}<I_t^s}}(X_s-I_t^s)\delta 
_{H_s}(dr).
\]
The process $\eta$ is the dual process of $\rho$ under $\N$ (see
Corollary 3.1.6 in \cite{dlg:rtlpsbp}). 

The next Lemma on time reversibility can easily be deduced
from Corollary 3.1.6 of  \cite{dlg:rtlpsbp} and the construction of
$m$. 

\begin{lem}
   \label{lem:reversib}
   Under  $\N$, the processes  $((\rho_s, \eta_s,
   \ind_{\{ m  _s=0\}}), s \in [0,  \sigma])$ and $((\eta_{(\sigma-s)-},
   \rho_{(\sigma-s)-}, $  $\ind_{\{ m  _ {(\sigma-s)-}=0\}}), s  \in [0,
   \sigma])$ have the same distribution.
\end{lem}

We present a Poisson  representation of $(\rho,\eta,m)$ under $\N$.  Let
$\mathcal{N}_0(dx\,   d\ell\,  du)$,   $\cn_1(dx\,   d\ell\,  du)$   and
$\cn_2(dx)$  be  independent  Poisson  point  measures  respectively  on
$[0,+\infty)^3$,  $[0,+\infty)^3$   and  $[0,+\infty)$  with  respective
intensity
\[
dx\,\ell\pi_0(d\ell)\, 
\ind_{[0,1]}(u)du, \quad dx\,\ell\pi_1(d\ell)\, 
\ind_{[0,1]}(u)du \quad\text{and}\quad \alpha_1dx.
\]
For every $a>0$, let us denote by $\mathbb{M}_a$ the law of the pair
$(\mu,\nu, m^\text{nod},m^\text{ske})$ of  measures on $\R_+$ with
finite mass defined by: 
for any $f\in \cb_+(\R_+)$  
\begin{align}
\label{def:mu_a}
\langle \mu,f\rangle  & =\int\left(\mathcal{N}_0(dx\, d\ell\, 
du)+\mathcal{N}_1(dx\, d\ell\, 
du)\right)  \, \ind_{[0,a]}(x) u\ell f(x) + \beta \int_0^a f(r)\; dr,\\
\label{def:nu_a}
\langle \nu,f\rangle  & =\int\left(\mathcal{N}_0(dx\, d\ell\, 
du)+\mathcal{N}_1(dx\, d\ell\, 
du)\right)  \, \ind_{[0,a]}(x) (1-u)\ell f(x)+ \beta \int_0^a f(r)\; dr,\\
\label{def:m_a}
\langle m^\text{nod} ,f\rangle  & =\int\mathcal{N}_1(dx\,  d\ell\, du)
\, \ind_{[0,a]}(x) u \ell 
f(x)\quad\text{and} \quad \langle m^\text{ske} ,f\rangle  =
\int\mathcal{N}_2(dx) \, \ind_{[0,a]}(x)
f(x). 
\end{align}
\begin{rem}
   \label{rem:W}
   In  particular $\mu(dr)+\nu(dr)$  is defined  as $\ind_{[0,a]}(r)
   d\xi_r$,   where  $\xi$  is   a  subordinator  with   Laplace  exponent
   $\psi'-\alpha$.
\end{rem}

We    finally    set   $\mathbb{M}=\int_0^{+\infty}da\,    \expp{-\alpha
  a}\mathbb{M}_a$.  Using the  construction of the snake, it  is easy to
deduce  from  Proposition  3.1.3  in \cite{dlg:rtlpsbp},  the  following
Poisson representation.
\begin{prop}\label{prop:poisson_representation2}
For every non-negative measurable function $F$ on $\cm_f(\R_+)^4$, 
\[
\N\left[\int_0^\sigma      F(\rho_t,     \eta_t,     m_t)      \;     dt
\right]=\int\mathbb{M}(d\mu\, d\nu\,dm)F (\mu, \nu,m) ,
\]
where $m=(m^\text{nod},m^\text{ske})$ and $\sigma=\inf\{s>0; \rho_s=0\}$
denotes the length of the excursion.
\end{prop}

\section{The pruned exploration process}\label{sec:pruned}

We  define the  following continuous
additive functional  of the process $((\rho_t,m_t),t\ge  0)$: 
\begin{equation}
   \label{eq:def-At}
A_t=\int_0^t \ind_{\{m_s=0\}}\; ds\quad\text{for $t\geq
0$.}
\end{equation}
\begin{lem}
\label{lem:A_s=0}
   We have the following properties.
\begin{itemize}
   \item[(i)] For $\lambda>0$, $\N[1  -\expp{-\lambda
     A_\sigma}]={\psi_0}^{-1} (\lambda)$.
   \item[(ii)] $\N$-a.e. 0 and $\sigma$ are points of increase for
   $A$. More precisely, $\N$-a.e. for all $\varepsilon>0$, we have 
$A_\varepsilon>0$ and $A_\sigma-A_{(\sigma-\varepsilon) \vee 0}>0$.
   \item[(iii)] If $\lim_{\lambda \rightarrow \infty }
     \phi_1(\lambda)=+\infty $, then $\N$-a.e. the set $\{s; m_s\neq
     0\}$ is dense in 
   $[0,\sigma]$. 
\end{itemize}
\end{lem}

\begin{proof}
 We first prove (i). 
Let  $\lambda>0$. Before  computing $v=\N[1  -\exp{-\lambda A_\sigma}]$,
notice  that  $A_\sigma\leq  \sigma$  implies, thanks  to
\reff{eq:N_s}, that 
$v\leq \N[1 -\exp{-\lambda \sigma}] =\psi^{-1}(\lambda)<+\infty$.  We have
\[
   v= \lambda \N\left[\int_0^\sigma dA_t\, \expp{-\lambda \int_t^\sigma
       dA_u}\right] 
=\lambda \N\left[\int_0^\sigma dA_t\, \E^*_{\rho_t, 0}[ \expp{-\lambda
  A_\sigma}] \right],
\]
where  we  replaced $\expp{-\lambda  \int_t^\sigma  dA_u}$  in the  last
equality  by   $\E^*_{\rho_t,  m_t}[  \expp{-\lambda   A_\sigma}]$,  its
optional projection,  and used  that $dA_t$-a.e. $m_t=0$. In  order to
compute this  last expression, we  use the decomposition of  $\cs$ under
$\P^*_{\mu}$ according  to excursions  of the  total mass  of $\rho$
above  its  minimum, see  Lemma  \ref{lem:dlg-decomp}.   Using the  same
notations as in  this lemma, notice that under  $\P^*_{\mu}$, we have
$A_\sigma=A_\infty =\sum_{i\in  \ck} \bar A^i_\infty  $, where for  every $T\ge
0$,
\begin{equation}
   \label{eq:def-Ai}
\bar A^i_T =\int_0^T \ind_{\{\bar m^i_t=0\}} dt.
\end{equation}
By Lemma \ref{lem:dlg-decomp}, we get 
\[
\E^*_{\mu}[ \expp{-\lambda A_\sigma}] = \expp{ - \langle \mu,1\rangle  \N[ 1-
  \exp{-\lambda A_\sigma}]}=\expp{-v\langle \mu,1\rangle }.
\]
We have
\begin{align}
\label{eq:calcul-v}
v=\lambda  \N\Big[\int_0^\sigma dA_t \, \expp{ -v \langle
  \rho_t,1\rangle }\Big] 
&= \lambda \N\Big[\int_0^\sigma dt \,\ind_{\{m_t=0\}} \expp{ -v
    \langle \rho_t,1\rangle }\Big] \\
\nonumber
&=\lambda\int_0^{+\infty}da\, \expp{-\alpha   a}\mathbb{M}_a [\ind_{\{m=0\}} \expp{ -v
    \langle \mu,1\rangle }   ]\\
\nonumber
&=\lambda\int_0^{+\infty}da\, \expp{-\alpha   a  } \exp
\left\{
- \alpha_1 a - \int_0^a dx
\int_0^1   du \int _{(0,\infty )} \ell_1\pi_1(d\ell_1) 
\right\}\\
\nonumber
&\hspace{1cm} \exp\Big\{-\beta v a - \int_0^a dx
\int_0^1   du \int _{(0,\infty )} \ell_0\pi_0(d\ell_0) \left( 1- \expp{
   -v u\ell_0  }\right) \Big\} \\
\label{eq:calc-v0}
&=\lambda\int_0^{+\infty}da\, \exp\left\{- a   \int_0^1 du\; \psi'_0(uv)
\right\} \\ 
\label{eq:calc-v}
&=\lambda\frac{v}{\psi_0(v)}, 
\end{align}
where  we used  Proposition  \ref{prop:poisson_representation2} for  the
third   and  fourth  equalities,   and  for   the  last   equality  that
$\alpha_0=\alpha+ \alpha_1 +  \int _{(0,\infty )} \ell_1\pi_1(d\ell_1) $
as well as
\begin{equation}
   \label{eq:def-psi1}
\psi_0'(\lambda)= \alpha_0+ \int_{(0,\infty )} \pi_0(d\ell_0)\; 
\ell_0(1-\expp{-\lambda \ell_0 }).
\end{equation}
Notice   that    if   $v=0$,   then    \reff{eq:calc-v0}   implies   $v=
\lambda/\psi_0'(0)  $, which  is  absurd since  $\psi'_0(0)=\alpha_0>0$
thanks to  \reff{eq:alpha0>0}. 
Therefore we have $v\in  (0,\infty)$, and we can divide \reff{eq:calc-v}
by $v$ to get $\psi_0(v)=\lambda$. This proves (i).

Now, we prove (ii). If we let $\lambda $ goes to infinity in (i) and use
that $\lim_{r\rightarrow \infty }  \psi_0(r)=+\infty $, then we get that
$\N[A_\sigma>0]=+\infty $.   Notice that for $(\mu,\Pi)\in  \S$, we have
under $\P^*_{\mu,\Pi}$,  $A_\infty \geq \sum_{i\in  \ck} \bar A^i_\infty
$,   with  $\bar   A^i$   defined  by   \reff{eq:def-Ai}.   Thus   Lemma
\ref{lem:dlg-decomp}  implies  that  if  $\mu\neq 0$,  then  $\P^*_{\mu,
  \Pi}$-a.s.   $\ck $  is  infinite and  $A_\infty  >0$.  Using  the
Markov property at time $t$ of the snake under $\N$, we get that for any
$t>0$,  $\N$-a.e.  on  $\{\sigma>t\}$, we  have  $A_\sigma-A_t>0$.  This
implies that $\sigma$ is $\N$-a.e. a  point of increase of $A$.  By time
reversibility, see Lemma \ref{lem:reversib},  we also get that $\N$-a.e.
$0$ is a point of increase of $A$. This gives (ii).

If  $\alpha_1>0$  then the  snake  $((\rho_s,  W_s),  s\geq 0)$  is  non
trivial.  It is  well known  that, since  the L\'evy  process $X$  is of
infinite variation,  the set  $\{s; \exists t\in [0,  H_s),\ W_s(t)\neq 0\}$ is $\N$-a.e. dense in $[0,\sigma]$. This implies that
$\{s; m_s\neq 0\}$ is $\N$-a.e. dense in $[0,\sigma]$. 

If $\alpha_1=0$ and
$\pi_1((0,\infty ))=\infty  $, then the  set $\cj^1$ of jumping  time of
$X$ is $\N$-a.e. dense in $[0,\sigma]$. This also  implies that
$\{s; m_s\neq 0\}$ is $\N$-a.e. dense in $[0,\sigma]$. 

If $\alpha_1=0$ and
$\pi_1((0,\infty ))<\infty  $, then the  set $\cj^1$ of jumping  time of
$X$ is $\N$-a.e. finite. This implies that 
$\{s; m_s\neq 0\}\cap [0,\sigma]$ is $\N$-a.e. a finite union of
intervals, which, thanks to (i), is not dense in $[0,\sigma]$. 

To get (iii), notice that $\lim_{\lambda \rightarrow\infty }
\phi_1(\lambda)=\infty $ if and only if  $\alpha_1>0$ or
$\pi_1((0,\infty ))=\infty  $. 
\end{proof}

We set $C_t=\inf \{  r>0; A_r > t\}$ the  right continuous inverse of $A$, with
the  convention  that  $\inf\emptyset=\infty   $.   
  {F}rom    excursion    decomposition,     see    Lemma
\ref{lem:dlg-decomp}, (ii) of  Lemma \ref{lem:A_s=0} implies
the following Corollary.
\begin{cor}
\label{cor:C_0=0}
For any initial measures $(\mu,\Pi)\in \S$,   $\P_{\mu,\Pi}$-a.s. the process
$(C_t,t\geq 0)$ is finite.  If $m_0=0$, then
$\P_{\mu,\Pi}$-a.s.  $C_0=0$.
\end{cor}

  We   define  the   pruned  exploration   process   $\tilde  \rho=(\tilde
\rho_t=\rho_{C_t}, t\geq 0)$ and the pruned marked exploration process $\tilde
\cs= ( \tilde \rho, \tilde m)$, where $\tilde m=(m_{C_t}, t\geq 0)=0$.  
Notice  $C_t$  is  a
$\cf$-stopping  time for  any $t\geq  0$ and is finite a.s. from
Corollary \ref{cor:C_0=0}.
Notice  the process $\tilde \rho$, and thus
the  process $\tilde  \cs$,  is  c\`ad-l\`ag.  We  also  set $\tilde  H_t=
H_{C_t}$ and $\tilde \sigma=\inf \{t>0; \tilde \rho_t=0\}$. 

Let $\tilde \cf=(\tilde \cf_t, t\geq 0)$ be the filtration generated by the
pruned  marked exploration process $\tilde  \cs$ completed  the usual  way.  In
particular   $\tilde \cf_t\subset   \cf_{C_t}$,   where   if   $\tau$   is   an
$\cf$-stopping time, then $\cf_\tau$ is the $\sigma$-field associated with
$\tau$.

We     are     now     able     to     restate     precisely     Theorem
\ref{thm:law_pruned_bref}. Let  $\rho^{(0)}$ be the  exploration process
of a L\'evy process with Laplace exponent $\psi_0$.
\begin{theo}\label{thm:law_pruned}
For every finite measure $\mu$, the law of the pruned process $\tilde
\rho$ under $\P_{\mu,0}$ is the law of the exploration process
$\rho^{(0)}$ associated with a L\'evy process with Laplace exponent
$\psi_0$ under $\P_\mu$.
\end{theo}
The proof of this  Theorem is given
at the end of Section \ref{sec:law_pruned}. 



\section{A special Markov property}\label{sec:Markov_special}

In  order to  define the  excursions  of the  marked exploration process away
from 
$\{s\ge  0;\ m_s=0\}$,  we define  $O$ as  the interior  of  $\{s\ge 0,\
m_s\ne 0\}$. We shall see that the complementary of $O$ has positive Lebesgue
measure. 

\begin{lem}
\label{lem:O}
\begin{itemize}
\item[(i)] If the set $\{s\ge  0,\ m_s\ne 0\}$ is non empty then,
  $\N$-a.e.  $O$ is non empty.
\item[(ii)] If we have  $\displaystyle \lim_{\lambda\rightarrow \infty } \phi_1(\lambda)=\infty $,
then $\N$-a.e. the open set $O$ is dense in $[0,\sigma]$.
\end{itemize}
\end{lem}

\begin{proof}
  For any element  $s'$ in $\{s\ge 0,\ m_s\ne  0\}$, there exists $u\leq
  H_{s'}$ such that $m_{s'}([0,u])\ne 0$ and $\rho_{s'}([u,\infty ))>0$.
  Then we consider  $\tau_{s'}=\inf \{t>s',\rho_t([u,\infty ))=0\}$.  By
  the right  continuity of $\rho$, we have $\tau_{s'}>s'$  and the snake
  property \reff{eq:snake-prop} implies that $\N$-a.e.
  $(s',\tau_{s'})\subset O$.

Use (iii) of  Lemma \ref{lem:A_s=0} to get the last part. 
\end{proof}

We   write   $O=\bigcup_{i\in   \ci}(\alpha_i,\beta_i)$   and   say   that
$(\alpha_i,\beta_i)_{i\in \ci}$  are the excursions intervals  of the
marked exploration process $\cs=(\rho, m)$ away from $\{s\ge 0,\ m_s=0\}$.  For every
$i\in \ci$, let us  define the measure-valued process
$\cs^i=(\rho^i,m^i)$. For every $f\in \cb_+(\R_+)$,  $t\ge 0$, we set
\begin{equation}
\label{eq:def-ri-mi}
\begin{aligned}
   \langle \rho_t^i,f\rangle                                                             &
=\int_{[H_{\alpha_i},+\infty)}f(x-H_{\alpha_i})\rho_{(\alpha_i+t)\wedge
\beta_i}(dx)\\                        
\langle (m^\text{a})_t ^i,f\rangle                         &
=\int_{(H_{\alpha_i},+\infty)}f(x-H_{\alpha_i})m_{(\alpha_i+t)\wedge  
\beta_i}^\text{a} (dx) \quad\text{with a$\in
\{\text{nod},\text{ske}\}$} 
\end{aligned}
\end{equation}
and  $m^i_t=((m^\text{nod})^i_t,  (m^\text{ske})^i_t)$.  Notice
that  the  mass 
located at $H_{\alpha_i}$ is kept, if there is any, in the definition of
$\rho^i$  whereas  it  is  removed  in  the  definition  of  $m^i$.   In
particular,  if $\Delta_{\alpha_i}>0$,  then
$\rho^i_0=\Delta_{\alpha_i} \delta_ 
{0}$ and for every $t< \beta_i-\alpha_i$, the measure $\rho^i_t$ charges
0.   On   the  contrary,   as   $m^i_0=0$,   we   have,  for   every   $t<
\beta_i-\alpha_i$, $m^i_t(\{0\})=0$.

Let  $\tilde \cf_\infty  $ be  the $\sigma$-field  generated  by $\tilde
\cs=((\rho_{C_t},  m_{C_t}), t\geq 0)$. Recall that $\P_{\mu,  \Pi }^{*}(d\cs)$
denotes the  law of the marked exploration process  $\cs$ started at  $(\mu,\Pi)\in \S$ and stopped
when  $\rho$  reaches  0.    For  $\ell\in(0,+\infty)$,  we  will  write
$\P_\ell^*$ for $\P_{\ell\delta_0,0}^*$.

If  $Q$ is  a measure  on $\S$  and $\varphi$  is a  non-negative measurable
function defined on the measurable space $\R_+\times \S\times \S$, we
denote by
\[
Q[\varphi(u,\omega,\cdot)]=\int _{\S}\varphi(u,\omega,\cs)Q(d\cs).
\]
In other words, the integration concerns only the third component of the
function $\varphi$.

We can now state the Special Markov Property.
\begin{theo}[Special Markov property]
\label{th:SMP}
Let  $\varphi$   be  a  non-negative  measurable   function  defined  on
$\R_+\times \cm_f(\R_+)\times \S$.
Then, we have $\P$-a.s.
\begin{multline}
   \label{eq:MS}
\E\left[\exp\left(-\sum_{i\in
    \ci}\varphi(A_{\alpha_i},\rho_{\alpha_i-},\cs^i)\right)\biggm|\tilde \cf_\infty
\right]\\ 
=\exp\left(-\int_0^\infty 
    du\,\alpha_1
    \N\left[1-\expp{-\varphi(u,\mu,\cdot)}\right]_{|\mu=\tilde\rho_u} 
  \right) \\
\exp\left(-\int_0^\infty 
    du\,\int_{(0,\infty )} \pi_1(d\ell) \left(1-
     \E^*_\ell[\expp{-\varphi(u,\mu,\cdot)} ]_{|\mu=\tilde\rho_u}
   \right) \right).   
\end{multline}
In other words, the law under $\P$ of the  excursion  process
$\displaystyle \sum 
_{i\in \ci}\delta_{( 
A_{\alpha_i},\rho_{\alpha_i-},  \cs^i)}(du\, d\mu\, d\cs) $,  given
$\tilde \cf_\infty 
$,  is  the  law of  a  Poisson  point  measure  with  intensity 
$$
\ind_{ \{u\geq 0\}} du  \;
\delta_{ \tilde \rho_u} (d\mu)\; \left(\alpha_1 \N(d\cs) +\!
  \int_{(0,\infty )}\! \pi_1(d\ell) \P^*_\ell (d\cs) \right).$$
\end{theo}

Informally speaking, this Theorem gives the distribution of the marked
exploration process ``above'' the pruned CRT. 
The end of this section is now devoted to its  proof.

Let    us    first    remark    that,    if    $\lim_{\lambda\to+\infty}
\phi_1(\lambda)<+\infty$, we  have $\alpha_1=0$ and $\pi_1$  is a finite
measure.  Hence,  there is no  marks on the  skeleton and the  number of
marks on  the nodes  is finite  on every bounded  interval of  time. The
proof  of Theorem  \ref{th:SMP} in  that case  is easy  and left  to the
reader.    For   the   rest    of   this   Section,   we   assume   that
$\lim_{\lambda\to+\infty}\phi_1(\lambda)=+\infty$.

\subsection{Preliminaries}
Fix $t>0$  and $\eta>0$. For $\cs=(\cs_s=(\rho_s,m_s),  s\geq 0)$, 
we set  $B=\{\sigma(\cs)= +\infty \}\cup  \{T_\eta(\cs)=+\infty \}\cup\{
L_\eta(\cs)=-\infty \}$ where
$\sigma(\cs)=\inf\{s>0; \rho_s=0\}$, $ T_\eta(\cs) = \inf\{s\geq \eta ; \langle 
\rho_s,1 \rangle  \geq \eta 
 \} $  and $
L_\eta(\cs) = \sup\{s\in [0, 
\sigma(\cs)] ; \langle 
\eta_s,1 \rangle  \geq \eta   \}$,
with the convention $\inf \emptyset=+\infty $ and $\sup \emptyset= -\infty$. 

We consider non-negative
bounded  functions  $\varphi$  satisfying  the  assumptions  of  Theorem
\ref{th:SMP} and these four conditions:
\begin{itemize}
\item[($h_1$)] $\varphi(u,\mu,\cs)=0$  for any $u\geq t$. 
   \item[($h_2$)] $u\mapsto \varphi(u,\mu,\cs)$ is uniformly
   Lipschitz (with a constant that does not depend on $\mu$ and $\cs$).
   \item[($h_3$)] $\varphi(u,\mu,\cs)=0$ on  $B$; and 
     if $\cs \in B^c$ then
     $\varphi(u,\mu,\cs)$ depends on $\cs$ only through $(\cs_u,u\in
     [T_\eta , L_\eta])$. 
   \item[($h_4$)]   The  function  $\mu\mapsto   \varphi(u,\mu,\cs)$  is
     continuous  with  respect to  the  distance $D(\mu,\mu')+  |\langle
     \mu,1 \rangle -\langle \mu',1 \rangle|$ on $\cm_f(\R_+)$, where $D$
     is a distance  on $\cm_f(\R_+)$ which defines the  topology of weak
     convergence.
\end{itemize} 

\begin{lem}
  Let $\varphi$ satisfies $(h_1-h_3)$ and  let $w$
  be defined
  on $(0,\infty )\times [0,\infty )\times \cm_f(\R_+)$ by
\[
w(\ell,u,\mu)=     \E^*_\ell[\expp{-\varphi(u,\mu,\cdot)} ] .
\]
Then, for $\N-a.e.$ $\mu\in \cm_f(\R_+)$, the function $(\ell,u)\mapsto w(\ell,u,\mu)$ is uniformly continuous on $(0,\infty )\times [0,\infty )$. 
\end{lem}

\begin{proof}
Let $u>0$ and $\ell'>\ell$. If we set $\tau_\ell=\inf\{t\ge 0,\ \rho_t(\{0\})=\ell\}$ we have, by
the strong Markov property at time $\tau_\ell$ and assumption $(h_3)$, that 
$$\E_{\ell'}^*\left[\expp{-\varphi(u,\mu,\cdot)}\right]=\E_{\ell'}^*\left[\ind_{\{T_\eta>\tau_\ell\}}\E_\ell^*\left[\expp{-\varphi(u,\mu,\cdot)}\right]\right]+\E_{\ell'}^*\left[\expp{-\varphi(u,\mu,\cdot)}\ind_{\{T_\eta\le
    \tau_\ell\}}\right].$$
Therefore,
\begin{align*}
\bigl|w(\ell',u,\mu)-w(\ell,u,\mu)| & \le\E_{\ell'}^*\left[\ind_{\{T_\eta\le\tau_\ell\}}\E_\ell^*\left[\expp{-\varphi(u,\mu,\cdot)}\right]\right]+\E_{\ell'}^*\left[\expp{-\varphi(u,\mu,\cdot)}\ind_{\{T_\eta\le
    \tau_\ell\}}\right]\\
& \le 2\P_{\ell'}^*(T_\eta\le
    \tau_\ell)\\
&=2\P_{\ell'-\ell}^*(T_\eta<+\infty).
\end{align*}
Using Lemma
\ref{lem:dlg-decomp}, for $\ell'-\ell <\eta$, we get 
\[
|w(\ell',u,\mu)-w(\ell,u,\mu)|\le 2\left(1-\expp{-(\ell'-\ell) \N[T_\eta<\infty
    ]}\right) .  
\]
Since $\N[T_\eta<\infty ]<\infty $, we then deduce there exists a finite
constant $c_\eta$  s.t. for all function $\varphi$ satisfying
$(h_3)$, 
\begin{equation}
   \label{eq:c-eta}
|w(\ell', u,\mu)-w(\ell,u,\mu) |\leq  c_\eta |\ell'-\ell|. 
\end{equation}

The absolute continuity with respect to $u$ is a direct consequence of
assumptions $(h_1-h_2)$.
\end{proof}

\subsection{Stopping times}
\label{sec:stop}
Let  $R(dt,du)$ be  a  Poisson  point measure  on  $\R_+^2$ (defined  on
$(\S,\cf)$)  independent  of $\cf_\infty$  with  intensity the  Lebesgue
measure. We  denote by $\cg_t$ the $\sigma$-field  generated by $R(\cdot
\cap   [0,t]\times\R_+)$.   For   every  $\varepsilon>0$,   the  process
$R_t^\varepsilon:=R([0,t]\times [0,1/\varepsilon])$ is a Poisson process
with intensity $1/\varepsilon$.  We denote by $(e_k^\varepsilon,k\ge 1)$
the time intervals between  the jumps of $(R_t^\varepsilon,t\ge 0)$. The
random  variables $(e_k^\varepsilon,  k\geq 1)$  are  i.i.d. exponential
random  variables  with  mean  $\varepsilon$,  and  are  independent  of
$\cf_\infty$. They define  a mesh of $\R_+$ which is  finer and finer as
$\varepsilon$ decreases to 0.

{F}or
$\varepsilon>0$,   we    consider  
$T^\varepsilon_0=0$, $M^\varepsilon_0=0$  and for $k\geq 0$,
\begin{equation}
   \label{eq:def-MST}
\begin{aligned}
M_{k+1}^\varepsilon&=\inf\{i>M_k^\varepsilon; m_{T^\varepsilon_k+
  \sum_{j=M^\varepsilon_k +1 }^i  e_j^\varepsilon}  \neq 0\}, \\ 
   S_{k+1}^\varepsilon & = T_{k}^\varepsilon+
\sum_{j=M^\varepsilon_k +1 }^{M_{k+1}^\varepsilon}  e_j^\varepsilon,\\
T_{k+1}^\varepsilon & =\inf\{ s>S_{k+1}^\varepsilon ; \;m_s=0   \},
\end{aligned}
\end{equation}
with  the  convention  $\inf  \emptyset  =+\infty  $. For every $t\ge
 0$, we set $\displaystyle \tau^\varepsilon_t=\int_0^t ds \;
 \ind_{\bigcup _{k\geq 1} [T^\varepsilon_k, S^\varepsilon_{k+1})}(s)$
 and 
\begin{equation}\label{eq:def-fe}
\cf_t^e=\sigma(\cf_t\cup \cg_{\tau_t^\varepsilon}).
\end{equation}
 Notice that
$T^\varepsilon_k$ and $S^\varepsilon_k$ are $\cf^e$-stopping times. 

Now  we
introduce  a notation for  the process  defined above  the marks  on the
intervals  $\left[ S_k^\varepsilon ,  T_k^\varepsilon \right]$.  We set,
for $a\geq 0$, $\bar{H}_a$ the level of the first mark, $\rho_a^{-}$ the
restriction of $\rho_a$ strictly below it and $\rho_a^+$ the restriction
of $\rho_a$ above it: 
\begin{equation}
   \label{eq:def-Htilde}
\bar{H}_a=\sup\{ t>0 , m_a([0,t])=0 \}, \quad \rho_a^- =\rho_a(\cdot\cap[0,
  \bar H_a ))
\end{equation}
and $\rho_a^+$ is defined by $\rho_a=[\rho_a^-, \rho_a^+]$,  that is for
any $f\in \cb_+(\R_+)$,  
\begin{equation}
   \label{eq:rho+}
   \langle\rho^{+}_a,
   f  \rangle=\int _{[\bar  H_a, \infty  )} f(r-  \bar  H_a)\; \rho_a(dr).
\end{equation}
For $k\geq 1$
and  $\varepsilon >0$ fixed,  we define  $\cs^{k,\varepsilon }  = \left(
  \rho^{k,\varepsilon }  , m^{k,\varepsilon }  \right)$ in the
  following way: for  $s >0$
and $f\in \cb_+(\R_+)$
\begin{align*}
 \rho_s^{k,\varepsilon }&= \rho^+_{(S_k^\varepsilon +s) \wedge
  T_k^\varepsilon },\\
\langle  (m^\text{a})_s^{k,\varepsilon },f \rangle
&=\int_{(\bar {H}_{S_k^\varepsilon },+\infty)}
f(r-\bar {H}_{S_k^\varepsilon }) m^\text{a}_{(S_k^\varepsilon +s)\wedge
  T_k^\varepsilon }(dr),  \quad\text{with a$\in
\{\text{nod},\text{ske}\}$},
\end{align*}
and $m_s^{k,\varepsilon }=((m^\text{nod})_s^{k,\varepsilon } ,
(m^\text{ske})_s^{k,\varepsilon }) $.
Notice that $\rho_s^{k,\varepsilon }(\{0\})=\rho_{S_k^\varepsilon
  }(\{ \bar {H}_{S_k^\varepsilon }\})$. 
For $k\geq 1$, we  consider the $\sigma$-field
$\cf^{(\varepsilon), k}$ generated by the family of processes
$\displaystyle \left(\cs_{(T_\ell^\varepsilon+s)\wedge
S_{\ell+1}^\varepsilon   -},\   s>0\right)_{ \ell\in  \{0, \ldots,
k-1\}}$.

Notice that 
for $k\in \N^*$
\begin{equation}
   \label{eq:inclusion}
\cf^{(\varepsilon),k} \subset
\cf^e_{S^\varepsilon_k}. 
\end{equation}

\subsection{Approximation of the functional}

\begin{figure}[H]
\psfrag{alphai}{$a_i$}
\psfrag{betai}{$b_{i}$}
\psfrag{alphai0}{$a_{i_0}$}
\psfrag{betai0}{$b_{i_0}$}
\psfrag{Ht}{$H_t$}
\begin{center}
\includegraphics[width=13cm]{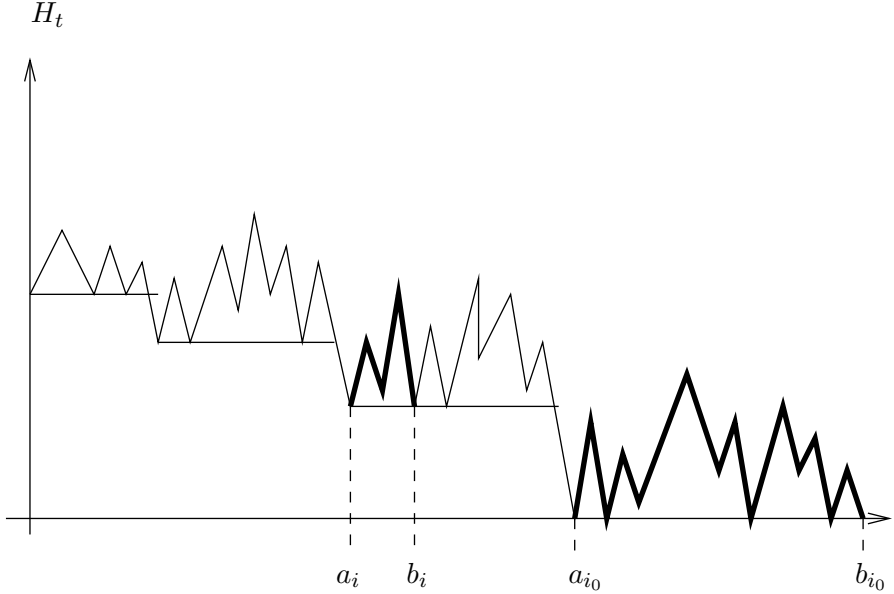}
\caption{Definition of the different excursions\label{fig:exc}}
\end{center}
\end{figure}

Let $\cs$ be a marked exploration  process and $g$ be a function defined
on $\S$.  We decompose the path of $\rho$ according to the excursions of
the   total  mass   of  $\rho$   above   its  minimum   as  in   Section 
\ref{sec:poisson}, with a slight difference if the initial measure $\mu$
charges $\{0\}$. More precisely, we perform the same decomposition as in
Section\ref{sec:poisson}  until   the  height  process   reaches  0.  If
$\mu(\{0\})=0$, then $H_t=0\iff \rho_t=0$ and the decompositions are the
same.  If  not,  there remains  a  part  of  the  process which  is  not
decomposed and is gathered in  a single excursion (defined as $(a_{i_0},
b_{i_0})$  in Figure  \ref{fig:exc}).  We set  $Y_t=  \langle \rho_t,  1
\rangle$ and  $J_t=\inf_{0\leq u\leq t}  Y_t$. Recall that  $(Y_t, t\geq
0)$ is  distributed under  $\P_\mu^*$ as a  L\'evy process  with Laplace
exponent $\psi$  started at $\langle  \mu,1 \rangle$ and killed  when it
reaches $0$. Let $(a_i, b_i)$, $i\in \ck$, be the intervals excursion of
$Y-J$   away    from   $0$.   For   every   $i\in    \ck$,   we   define
$h_i=H_{a_i}=H_{b_i}$.  We set $\tilde\ck=\{i\in  \ck; h_i>0\}$  and for
$i\in \tilde\ck$  let $\bar\cs^i=(\bar\rho^i,  \bar m^i)$ be  defined by
\reff{eq:rho-min}  and \reff{eq:m-min}.  If the initial measure $\mu$
does not charge $0$, we have $\tilde \ck=\ck$ and we set $\ck^*=\tilde
\ck=\ck$. If the initial measure $\mu$ charges $0$, we consider
$i_0\not\in  \tilde\ck$ and 
set   $\ck^*=\tilde\ck\cup\{i_0\}$,    $a_{i_0}=\inf\{a_i;   i\in   \ck,
h_i=0\}$,     $b_{i_0}=\sup\{b_i;     i\in     \ck,     h_i=0\}$     and
$\bar\cs^{i_0}=(\bar\rho^{i_0} ,\bar m^{i_0})$ with
\begin{align*}
\langle \bar\rho_t^{i_0},f\rangle
& =\int_{[0,+\infty)}f(x)\rho_{(a_{i_0}+t)\wedge
\tilde b_{i_0}}(dx)\\                        
\langle (\bar m^\text{a})_t^{i_0},f\rangle                         &
=\int_{(0,+\infty)}f(x)m^\text{a} _{(a_{i_0}+t)\wedge
b_{i_0}}(dx)\quad\text{with a$\in
\{\text{nod},\text{ske}\}$} ,
\end{align*}
and     $\bar     m^{i_0}_t=((\bar    m^\text{nod})^{i_0}_t,     (\bar
m^\text{ske})^{i_0}_t)$.  See Figure  \ref{fig:exc} to get the picture
of the different excursions.

We define 
\begin{equation}
   \label{eq:def-g*}
g^*(\cs)=\sum_{i\in \ck^*} g(\bar\cs^i). 
\end{equation}

\begin{lem}\label{lem:approx}
$\P$-a.s., we have, for $\varepsilon>0$ small enough,
\begin{equation}
   \label{eq:=cs}
\sum_{i\in
  \ci}\varphi(A_{\alpha_i},\rho_{\alpha_i-},\cs^i)=\sum_{k=1}^\infty
\varphi(A_{S_k^\varepsilon},\rho_{S_k^\varepsilon}^-,\cs^{k,\varepsilon})
=\sum_{k=1}^\infty 
  \varphi^*(A_{S_k^\varepsilon},\rho_{S_k^\varepsilon}^-,\cs^{k,\varepsilon}),
\end{equation}
where the sums have a finite number of non zero terms. 
\end{lem}

\begin{proof}
{\em First equality}. 
By assumptions $(h_1)$ and $(h_3)$, as $\N[T_\eta<+\infty]<+\infty$,
the set
$$\ci'=\{i\in \ci, \varphi(A_{\alpha_i},\rho_{\alpha_i-},\cs^i)\ne 0\}$$
is finite. Therefore, for $\varepsilon$ small enough, for every $j\in
\ci'$, the mesh defined by \reff{eq:def-MST} intersects the interval
$(\alpha_j,\beta_j)$: in other words, 
there exists an integer $k_j$ such that
$S_{k_j}^\varepsilon\in(\alpha_j,\beta_j)$ (and that integer is unique).

Moreover, for every $j\in\ci'$,  we can choose $\varepsilon$ small enough
so  that  $S_{k_j}^\varepsilon<T_\eta(\rho^j)$,  which gives  that,  for
$\varepsilon$ small enough,
$$\varphi(A_{\alpha_j},\rho_{\alpha_j-},\cs^j)=\varphi(A_{\alpha_j},\rho_{\alpha_j-},\cs^{k_j,\varepsilon}).$$
Finally, as the mark at $\alpha_j$ is still present at time
$S_{k_j}^\varepsilon$, the additive functional $A$ is constant on that
time interval and
$\rho_{\alpha_j-}=\rho_{S_{k_j}^\varepsilon}^-$. Therefore, we get
$$\varphi(A_{\alpha_j},\rho_{\alpha_j-},\cs^j)=\varphi(A_{S_{k_j}^\varepsilon},\rho_{S_{k_j}^\varepsilon}^-,\cs^{k_j,\varepsilon}).$$

{\em Second equality}.  Let $j\in\ci'$. We consider the decomposition of
$\cs^{k_j,\varepsilon}$ according  to $\rho^{k_j,\varepsilon}$ above its
minimum described at the beginning of this Subsection.  We must consider
two cases :

\medskip
First case : The mass at $\alpha_j$ is on a node.
Then, for $\varepsilon>0$ small enough, we have $T_\eta>a_{i_0}$
and
$$\varphi(A_{S_{k_j}^\varepsilon},\rho_{S_{k_j}^\varepsilon}^-,\cs^{k_j,\varepsilon})=\varphi(A_{S_{k_j}^\varepsilon},\rho_{S_{k_j}^\varepsilon}^-,\cs^{i_0})=\varphi^*(A_{S_{k_j}^\varepsilon},\rho_{S_{k_j}^\varepsilon}^-,\cs^{k_j,\varepsilon})$$
as all the terms in the sum that defines $\varphi^*$ are zero but for
$i=i_0$.

\medskip
Second case : The mass at $\alpha_j$ is on the skeleton.
In that case, we again can choose $\varepsilon$ small enough so that
the interval $[T_\eta,L_\eta]$ is included in one excursion interval
above the minimum of the exploration total mass process of $S^{k_j,
  \varepsilon}$. We then conclude as in the previous case.
\end{proof}


\subsection{Computation of the conditional expectation of the
  approximation}

\begin{lem}
\label{lem:ccea}
For every $\tilde \cf_\infty$-measurable non-negative random variable
$Z$, we have 
\[
   \E\left[ Z \exp\left( -\sum_{k=1}^\infty 
\varphi^*\left
      (\AAke ,\rho_{S_k^\varepsilon}^-,\cs^{k,\varepsilon }
    \right) \right) 
\right] 
=\E\left[Z\prod_{k=1}^\infty 
K_\varepsilon (\AAke,  
    \rho_{S_k^\varepsilon}^- )\right],
\]
where $\displaystyle \gamma =
\psi^{-1}\left(1/{\varepsilon }\right)$ and 
\begin{equation}
   \label{eq:Ke}
K_\varepsilon(r,  
    \mu )= \frac{\psi(\gamma)}{\phi_1(\gamma)}\frac{\gamma - v(r,\mu)
    }{ \psi(\gamma)- \psi(v(r,\mu)) } \left( \alpha_1 +   
  \int_0^1 du    \int_{(0, \infty )} \ell_1  \pi_1(d\ell_1)\;
  w(u\ell_1,r,\mu ) \expp{-\gamma(1-u)\ell_1} \right).
\end{equation}
where 
\begin{equation}
   \label{eq:wv}
w(\ell,r,\mu)=\E^*_{\ell } \left[ \expp{-\varphi( r
    ,\mu,\cdot)}\right] \quad\text{and}\quad 
v(r,\mu)= \N \left[1-\expp{-\varphi( r  ,\mu,\cdot)}  \right]. 
\end{equation}
\end{lem}

\begin{proof}
This proof is rather long and technical. We decompose it in several steps.

\bigskip
{\sl Step 1.} We introduce first a special form of the random
variable $Z$.

  Let  $p\ge 1$. Recall that $H_{t,t'}$ denotes the
  minimum of $H$ between $t$ and $t'$ and that $\bar H_t$ defined by
  \reff{eq:def-Htilde} represents the height of the lowest mark. We set 
\begin{align*}
\xi_d^{p-1}  &=\sup\left\{t>T_{p-1}^\varepsilon;\
H_t=H_{T_{p-1}^\varepsilon,S_p^\varepsilon}\right\},\\
\xi_g^p  &=\inf\left\{t>T_{p-1}^\varepsilon;\
H_t=\bar H_{S_p^\varepsilon}\quad\text{and}\quad H_{t,
  S^\varepsilon_p}=H_t \right\}.   
\end{align*}
$\xi_d^{p-1}$ is the time at which the height process
reaches its minimum over $[T^\varepsilon_{p-1}, S^\varepsilon_p]$. By
definition of $T_{p-1}^\varepsilon$, $m_{T_{p-1}^\varepsilon}=0$
(there is no mark on the linage of the corresponding individual). On
the contrary, $m_{S_p^\varepsilon}\ne 0$, $m_{S_p^\varepsilon}(\{\bar
H_{S_p^\varepsilon}\})\ne 0$ but $m_{S_p^\varepsilon}([0, \bar
H_{S_p^\varepsilon}))=0$. In other words, at time $S_p^\varepsilon$,
some mark exists and the lowest mark is situated at height $\bar
H_{S_p^\varepsilon}.$ Roughly speaking,
$\xi_g^p$ is the time at which this lowest mark appears,
see figure \ref{fig:xi} to help understanding. Let us recall that, by
the snake property \reff{eq:snake-prop}, $m_{\xi_d^{p-1}}=0$ and
consequently, $\xi_d^{p-1}<\xi_g^p$ a.s.

\begin{figure}[H]
\psfrag{H}{$H$}
\psfrag{T}{$T_{p-1}^\varepsilon$}
\psfrag{S}{$S_{p}^\varepsilon$}
\psfrag{xid{p-1}}{$\xi^{p-1}_d$}
\psfrag{xig{p}}{$\xi^{p}_g$}
\psfrag{H}{$H$}
\psfrag{t}{$t$}
\psfrag{p}{}
\psfrag{p-1}{}
\psfrag{Hbar}{$\bar H_{S_p^\varepsilon}$}
\begin{center}
\includegraphics[width=13cm]{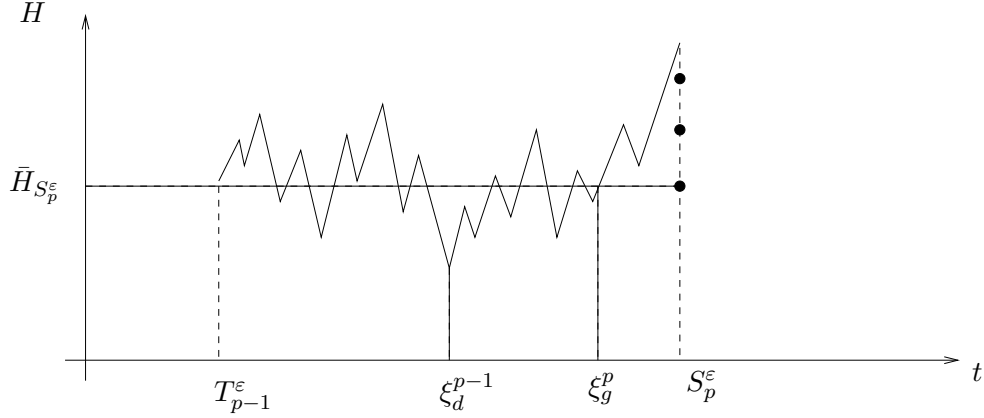}
\label{fig:xi}
\caption{Position of various random times}
\end{center}
\end{figure}

We consider a bounded non-negative random variable
  $Z$ of  the form $Z=Z_0 Z_1 Z_2 Z_3$, where $Z_0\in 
  \cf^{(\varepsilon),  p-1} $, 
$Z_1\in \sigma(\cs_u, T_{p-1}^\varepsilon\le u<\xi_d^{p-1})$,
$Z_2\in \sigma(\cs_u, \xi_d^{p-1}\le u<\xi_g^p)$
and
  $Z_3\in \sigma (\cs_{(T_k^\varepsilon+s)\wedge S_{k+1}^\varepsilon
    -},\ s\ge0,\ k\geq p)$ are bounded non-negative.

  \bigskip {\sl Step 2.} We apply  the strong Markov property to get rid
  of  terms which involve  $S^\varepsilon_p$ and  $T^\varepsilon_p$.  We
  first apply  the strong Markov property at  time $T_p^\varepsilon$ by
  conditioning with respect to $\cf^e_{T^\varepsilon_p}$. We
  obtain
\begin{multline*}
\E\left[Z\exp\left(-\sum_{k=1}^p
  \varphi^*(\AAke,\rho_{S_k^\varepsilon}^-,\cs^{k,\varepsilon}) 
  \right)\right]\\
=\E\left[Z_0Z_1Z_2\exp\left(-\sum_{k=1}^p
  \varphi^*(\AAke, \rho_{S_k^\varepsilon}^-,\cs^{k,\varepsilon})\right)
  \E_{\rho_{T_p^\varepsilon},0}
\left[Z_3\right]\right].
\end{multline*}
Recall notation \reff{eq:def-Htilde} and \reff{eq:rho+}. 
Notice    that
$\rho_{T_p^\varepsilon}=\rho_{S_p^\varepsilon}^-$,   and
consequently  $\rho_{T_p^\varepsilon}$  is  measurable with  respect  to
$\mathcal{F}_{S_p^\varepsilon}^e$.   So,  when we  use  the strong  Markov
property at time $S_p^\varepsilon$, we get thanks to \reff{eq:inclusion}
\begin{multline*}
\E\left[Z\exp\left(-\sum_{k=1}^p\varphi^*(A_{S_k^\varepsilon}
    ,\rho_{S_k^\varepsilon}^-,  
    \cs^{k,\varepsilon})\right)\right]\\  
=\E\left[Z_0Z_1Z_2\exp\left(-\sum_{k=1}^{p-1}\varphi^*(A_{S_k^\varepsilon}
    ,\rho_{S_k^\varepsilon}^-,  
    \cs^{k,\varepsilon})\right)
 \E_{\rho^{+}_{S^\varepsilon_p}    } ^*\left[\expp{-\varphi^*(
     A_{S_p^\varepsilon} ,\rho_{S_p^\varepsilon}^-,\cdot)}\right] 
   \E_{\rho_{S^\varepsilon_p}^{-}}[Z_3]\right].
\end{multline*}
Using the strong Markov property at time $T^\varepsilon_{p-1}$, and the
lack of memory for the exponential r.v.,  we get
\begin{multline}
\label{eq:N-EFG}
\E\left[Z\exp\left(-\sum_{k=1}^p\varphi^*(A_{S_k^\varepsilon}
    ,\rho_{S_k^\varepsilon}^-,  
    \cs^{k,\varepsilon})\right)\right]\\  
=\E\left[Z_0\exp\left(-\sum_{k=1}^{p-1}\varphi^*(A_{S_k^\varepsilon}
    ,\rho_{S_k^\varepsilon}^-,  
    \cs^{k,\varepsilon})\right)\phi\left(A_{S_{p-1}^\varepsilon},\rho_{S_{p-1}^\varepsilon}^-,\rho_{T_{p-1}^\varepsilon}\right)
 \right], 
\end{multline}
with
\begin{equation}
   \label{eq:def-phi-bmn}
\phi(b,\mu,\nu)=\E_\nu^*\left[Z_1Z_2\E_{\rho_{\tau'}^-}^*
  \left[\expp{-\varphi^*(b+A_{\tau'},\mu,\cdot)}\right]
  \E_{\rho_{\tau'}^-}^*[Z_3]\right], 
\end{equation}
where $\tau'$ is distributed under $\P_\nu^*$ as $S_1^\varepsilon$.

\bigskip
{\sl Step 3.}
We  compute the function $\phi$ given by \reff{eq:def-phi-bmn}. To simplify the formulas,
 we set
$$F(b',\mu')=\E_{\mu'}^*  \left[\expp{-\varphi^*(  b+b'
    ,\mu,\cdot)}\right],\qquad  G(\mu')= 
    \E_{\mu'}^*[Z_3]$$
(the dependence on $b$ and $\mu$ of $F$ is omitted) so that
\begin{equation}
   \label{eq:phi-bnm2}
\phi(b,\mu,\nu)=\E_\nu^*\left[Z_1Z_2F(A_{\tau'},
  \rho_{\tau'}^+)G(\rho_{\tau'}^-)\right]. 
\end{equation}

The proof of the following technical Lemma is postponed to the end of
this Sub-section. 
\begin{lem}
   \label{lem:lem-techn}
We set  $q(du, d\ell_1)=\alpha_1 \delta_{(0,0)}(dud\ell_1)+  du\; \ell_1
\pi_1(d\ell_1)$ and by convention $\pi(\{0\})=0$. We have:
\begin{equation}
   \label{eq:EGGFG1}
  \phi(b,\mu,\nu)
=  \E_{\nu}\left[Z_1 Z_2\frac{\Gamma_ F(A_{\tau'})}{ \Gamma_1}G(\rho_{\tau'}^-)
\right], 
\end{equation}
where for a non-negative function $f$  defined on $[0,\infty )\times\cm_f(\R_+)$
\[
\Gamma_f(a)
= \int_{[0,1]\times [0,\infty )} 
q(du, d\ell_1)\;  \int \M( d\rho', d\eta', dm')
\; \expp{-\gamma \langle \rho', 1 \rangle -
  \gamma u\ell_1}f(a,\eta'+(1-u)\ell_1\delta_0)
\]
and for $f=1$, $\Gamma_1$ does not depend on $a$.
\end{lem}

We now use the particular form of $F$ to compute $\Gamma_F$. 
Using \reff{eq:def-g*} and Lemma \ref{lem:dlg-decomp}, we get
\[
F(a,\mu')
=\E_{\mu'}^*  \left[\expp{-\varphi^*(b+a
     ,\mu,\cdot)}\right] \\
=\E^*_{\mu'(\{0\})} \left[ \expp{-\varphi(b+a
      ,\mu,\cdot)}\right] \expp{
-\mu'((0,\infty )) \N \left[1-\expp{-\varphi(b+a
      ,\mu,\cdot)}  \right]} .
\]
Using $w$ and $v$ defined in \reff{eq:wv}, we get 
\begin{multline*}
  \M_{s}\left[\expp{-\gamma \langle \rho,1
    \rangle-\gamma u\ell_1} 
F(a, \eta+ (1-u)\ell_1\delta_0)\right]\\
 \begin{aligned}
& =w((1-u)\ell_1,b+a,\mu) \expp{-\gamma u\ell_1} \M_{s}\left[\expp{-\gamma
    \langle \rho,1 
    \rangle} 
 \expp{
-v(b+a,\mu) \langle \eta, 1 \rangle } \right] \\
& =w((1-u)\ell_1, b+a,\mu ) \expp{-\gamma u\ell_1} \expp{ - s
  \left(\frac{\psi(\gamma)-\psi(v(b+a,\mu)) }{\gamma- v(b+a,\mu) } -\alpha \right)}.
\end{aligned}
\end{multline*}
We deduce that 
\begin{align*}
   \Gamma_{F}(a)
&=\frac{\gamma-v(b+a,\mu) }{ \psi(\gamma)-\psi(v(b+a,\mu)) } \left( \alpha_1 + 
  \int_0^1 du    \int_{(0, \infty )} \ell_1  \pi_1(d\ell_1)\;
  w(u\ell_1, b+a,\mu ) \expp{-\gamma(1-u)\ell_1} \right),
\end{align*}
and with $F=1$, $\displaystyle  \Gamma_1=\frac{\gamma}{
  \psi(\gamma)}  \frac{\phi_1(\gamma)}{\gamma}=\frac{\phi_1(\gamma)}{ 
  \psi(\gamma)} $.

Finally, plugging this formula in \reff{eq:EGGFG1} and using the
function $K_\varepsilon$ introduced in \reff{eq:Ke}, we have
\begin{equation}
\label{eq:phi3}
\phi(b,\mu,\nu)=\E_\nu[Z_1Z_2K_\varepsilon(b+A_{\tau'},\mu)G(\rho_{\tau'}^-)].
\end{equation}

\bigskip
{\sl Step 4.} Induction.

Plugging the expression \reff{eq:phi3} for $\phi$ in \reff{eq:N-EFG}, and using the arguments backward
from \reff{eq:N-EFG} we get 
\[
\E\left[Z\exp\left(-\sum_{k=1}^p\varphi^*(\AAke,
    \rho_{S_k^\varepsilon}^-,     \cs^{k,\varepsilon})\right)\right]
=
\E\left[Z\exp\left(-\sum_{k=1}^{p-1}\varphi^*(\AAke, 
    \rho_{S_k^\varepsilon}^-,     \cs^{k,\varepsilon})\right) K_\varepsilon
(A_{S_p^\varepsilon},  
    \rho_{S_p^\varepsilon}^- )\right]  .
\]
In particular, from  monotone  class Theorem,  this  equality  holds  for  any
non-negative $Z$ measurable  w.r.t.  the $\sigma$-field $\bar
\cf_\infty^\varepsilon 
=         \sigma((\cs_{C_t},         t\in         [A_{T^\varepsilon_{k}},
A_{S^\varepsilon_{k+1}}]),   k\geq   0)$.   Notice  that   $K_\varepsilon  (A_{S_p^\varepsilon}, \rho_{S_p^\varepsilon}^-   )$  is
measurable w.r.t.  $\bar \cf_\infty $.  So, we may iterate  the previous
argument  and let  $p$ goes  to  infinity to  finally get  that for  any
non-negative random variable $Z \in \bar  \cf_{\infty} $, we have
\[
\E\left[Z\exp\left(-\sum_{k=1}^\infty 
\varphi^*(A_{S_k^\varepsilon}
    ,\rho_{S_k^\varepsilon}^-,  
    \cs^{k,\varepsilon})\right)\right]  
=\E\left[Z\prod_{k=1}^\infty 
K_\varepsilon (\AAke,  
    \rho_{S_k^\varepsilon}^- )\right].
\]
Intuitively,  $\bar  \cf_{\infty}^\varepsilon  $ is  the  $\sigma$-field
generated by  $\tilde \cf_\infty $ and  the mesh $([A_{T^\varepsilon_k},
A_{S^\varepsilon_{k+1}}], k\geq 0)$. As $\bar \cf_{\infty}^\varepsilon $
contains $\tilde \cf_\infty $, the Lemma is proved.
\end{proof}

\begin{proof}[Proof of Lemma \ref{lem:lem-techn}]
We consider the Poisson decomposition of $\cs$ under $\P^*
_\nu $ given in Lemma \ref{lem:dlg-decomp}. 
Notice there exists a unique excursion $i_1\in \ck$
s.t. $a_{i_1}<\tau'<b_{i_1}$. 

By hypothesis on $Z_1$, under $\P^*_\nu$, we can write
$Z_1=\ch_1(\nu , \sum_{i\in \ck;
  a_i<a_{i_1}} \delta_{ h_i,
  \bar\cs^i} )$ for a measurable function $\ch_1$. We can also write
$Z_2=\ch_2(\rho_u,\xi_d^{0}\le u<\xi_g^{1})$ for a measurable
function $\ch_2$ as $m_u=0$ for $u\in [\xi_d^{0},\xi_g^1)$. 
Then, using compensation  formula  in  excursion  theory,  see
Corollary IV.11  in  \cite{b:pl}, we get 
\begin{equation}
   \label{eq:hGF}
    \phi(b,\mu,\nu)
=\E\left[\int  \nu(dv)\; \ind_{\displaystyle \{\tau'>\sum_{s<v}  \sigma(\cs^s) \}}\ch_1\left(\nu,\sum_{s<v}
\delta_{(s,\cs^s)}\right)h_{F}^{\nu} \left(r,\sum_{s<v} A_{\sigma(\cs^s)}(\cs^s) \right)
\right], 
\end{equation}
where  $\sum_{s} \delta_{s,\cs^s}$ is a Poisson point measure
with intensity $\nu(ds) \N[d\cs]$, $\sigma(\cs)=\inf\{r>0, \cs_r=0\}$,
$A_t(\cs^s)=\int_0^t dv'\; \ind_{\{m_{v'}(\cs^s)=0\}}$  and 
\[
h^{\nu}_{F}(r,a)=\N\left[
  F(a+ A_{\tau'},\rho^{+}_{\tau'}) 
    G([k_{r}{\nu},\rho_{\tau'}^-])\ch_2([k_r\nu,\rho_t], 0\le t<
    \xi_g^1)\ind_{\{\tau'<\sigma\}} \right].
\]

Let $(R_k, k\ge 1)$ be the increasing sequence of the jumping times of a
Poisson process of intensity $1/\varepsilon$, independent of
$\cs$. Then, by time-reversal, we have 
\begin{multline*}
h^{\nu}_{F}(r,a)=\N\Big[\sum_{k=1}^{+\infty}\ind_{\{m_{R_k}\ne
    0\}}\ind_{\{\forall k'>k,\ m_{R_{k'}}=0\}}F( a+
  A_\sigma -A_{R_k}  ,\eta^{+}_{R_k})\\ 
    G([k_{r}{\nu},\eta_{R_k}^-])\ch_2([k_r\nu,\eta_u],\tau_k<
    u\le\sigma)\Big],
\end{multline*}
where $\tau_k=\inf\{t>R_k;m_t=0\}$.
We then apply the strong Markov property at time $R_k$ and the Poisson
representation of the marked exploration process  to get 
\begin{multline*}
h^{\nu}_{F}(r,a)
=\N\Big[\sum_{k=1}^{+\infty}\ind_{\{m_{R_k}\ne
    0\}}
    G([k_{r}{\nu},\eta_{R_k}^-])\\
\E_{(\rho_{R_k}, m_{R_k})}^*\left[\ind_{\{\forall
    k'>0,\ m_{R_{k'}}=0\}}F(a+A_\sigma,\eta')  \ch_2([k_r\nu,\eta_u],\tau_0<
      u\le\sigma)\right]_{|\eta'=\eta^{+}_{R_k}}\Big],
\end{multline*}
where $\tau_0=\inf\{t>0;m_t=0\}$.   Now, let us remark  that, if $m_0\ne
0$, then $m_s\ne 0$ for $s\in [0,\tau_0]$ and $A_{\tau_0}=0$. Therefore,
$ m_{R_{1}}=0$ implies $R_1>\tau_0.$  The strong Markov property at time
$\tau_0$ gives, with $\eta'=\eta^{+}_{R_k}$, 
\begin{multline*}
   \ind_{\{m_{R_k}\neq 0\}} \E_{(\rho_{R_k}, m_{R_k})}^*  \left[\ind_{\{\forall
    k'>0,\ m_{R_{k'}}=0\}}F(a+A_\sigma, \eta' ) \ch_2([k_r\nu,\eta_u],\tau_0<
      u\le\sigma)\right]\\
 =\ind_{\{m_{R_k}\neq 0\}}
\P_{\rho_{R_k}^+}^*(R_1>\sigma)\E_{\rho_{R_k}^-}^*\left[\ind_{\{\forall 
    k'>0,\ m_{R_{k'}}=0\}}F(a+A_\sigma,\eta' ) \ch_2([k_r\nu,\eta_u],0<u\le
    \sigma)\right].
\end{multline*}
We have, using the Poisson representation of Lemma \ref{lem:dlg-decomp}
and \reff{eq:N_s}, that 
$$\P_{\rho_{R_k}^+}^*(R_1>\sigma)
=\E_{\rho_{R_k}^+}^*\left[\expp{-\sigma/\varepsilon} 
\right]= \expp{-\gamma\langle
  \rho_{R_k}^+,1\rangle},$$
as $\gamma=\psi^{-1}(1/\varepsilon)$. We
obtain 
\[
h_{F}^{\nu}(r,a)  =\N\left[\sum_{k=1}^{+\infty}\ind_{\{m_{R_k}\ne
    0\}}\widetilde  
    G(\rho_{R_k}^-, \eta_{R_k}^-, \rho_{R_k}^+, \eta_{R_k}^+)\right],
\]
where 
\begin{multline*}
  \widetilde G(\rho, \eta, \rho', \eta')= G([k_{r}{\nu},\eta])
\expp{-\gamma\langle
  \rho',1\rangle}\\
\E_{\rho}^*\left[\ind_{\{\forall 
    k'>0,\ m_{R_{k'}}=0\}}F( a+A_\sigma,\eta' ) \ch_2([k_r\nu,\eta],0<u\le
    \sigma)\right].
\end{multline*}
As $\sum_{k\geq 1} \delta_{R_k}$ is a Poisson point process with
intensity $1/\varepsilon$, we deduce that 
\[
h_{F}^{\nu}(r,a) =\frac{1}{\varepsilon}\N\left[\int_0^\sigma dt\,
  \ind_{\{m_{t}\ne 
    0\}} \widetilde  
    G(\rho_{t}^-, \eta_{t}^-, \rho_{t}^+, \eta_{t}^+)
  \right].
\]
Using Proposition \ref{prop:poisson_representation2},
we get
\[
 h_{F}^{\nu}(r,a)
=\inv{\varepsilon} \int_0^\infty da\; \expp{-\alpha a} \M_a
\left[ \ind_{\{m\ne 0\}}\widetilde
G(\mu^-,\nu^-,\mu^+,\nu^+)
\right].
\]

For $r>0$ and $\mu$ a measure on $\R_+$, let us define the measures
$\mu_{\ge r}$ and $\mu_{<r}$ by
\[
\langle\mu_{\ge r},f\rangle=\int f(x-r)\ind_{\{x\ge r\}}\mu(dx)\quad
\text{and}\quad \langle\mu_{< r},f\rangle=\int
f(x)\ind_{\{x<r\}}\mu(dx).
\]
Using Palm formula, we
get
\begin{multline*}
\M_a\left[\ind_{\{m\ne 0\}}\widetilde
G(\mu^-,\nu^-,\mu^+,\nu^+)  \right] \\
\hspace{-3cm}=\int_0^a dr\; 
\int_{[0,1]\times [0,\infty )}
q(du, d\ell_1)\; \\
\M_a\left[\ind_{\{m([0,r))=0\}}\widetilde
G(\mu_{<r},\nu_{<r},\mu_{\geq r}+ u\ell_1 \delta_0 ,\nu_{\geq r}+
(1-u)\ell_1 \delta_0 ) 
 \right]. 
\end{multline*}
Using the independence of the Poisson point measures, we get 
\begin{multline*}
\M_a\left[\ind_{\{m([0,r))=0\}}\widetilde
G(\mu_{<r},\nu_{<r},\mu_{\geq r}+ u\ell_1 \delta_0 ,\nu_{\geq r}+
(1-u)\ell_1 \delta_0 ) 
 \right]\\
=\int \M_r(d\mu,d\nu, dm)  \int \M_{a-r}( d\rho', d\eta', dm')
\ind_{\{m=0\}}\widetilde 
G(\mu,\nu,\rho'+ u\ell_1 \delta_0 ,\eta'+
(1-u)\ell_1 \delta_0 ) .
\end{multline*}
We deduce that 
\begin{multline*}
   h_{F}^{\nu}(r,a)
=\inv{\varepsilon} \int_{[0,1]\times [0,\infty )} 
q(du, d\ell_1)\; 
\int \M(d\rho,d\eta, dm)  \int \M( d\rho', d\eta',dm' )\\
\ind_{\{m=0\}}\widetilde 
G(\rho,\eta,\rho'+ u\ell_1 \delta_0 ,\eta'+
(1-u)\ell_1 \delta_0 ) .
\end{multline*}
Using this and \reff{eq:hGF} with similar arguments (in reverse order),
we obtain  \reff{eq:EGGFG1}. 
\end{proof}

\subsection{Computation of the limit}
  Recall
notation  of  Section  \ref{sec:stop}.   Let  $A^\varepsilon_s$  be  the
Lebesgue   measure    of   $[0,s]   \bigcap    \left(\bigcup_{k\geq   0}
  [T^\varepsilon_{k},S_{k+1}^\varepsilon]    \right)$.     The   process
$t\mapsto   \sup\{   i\in    \N;   \sum_{j=1}^i   e_j^\varepsilon   \leq
A^\varepsilon_t\}$ is a Poisson process with intensity $1/\varepsilon$ and
the process $s\mapsto N_{\varepsilon,t}$, where
\[
  N_{\varepsilon  ,t}=\sup\{k\in  \N ;  A_{S_k^\varepsilon
}^\varepsilon \leq  t  \}= 
\sup\{ k\in \N;  \sum_{j=1}^{M^\varepsilon_k}  e_j^\varepsilon \leq
A^\varepsilon_t\} ,
\]
is a marked Poisson process
with intensity $\P(m_\tau\neq 0)/\varepsilon$, where $\tau$ is an
 exponential random variable with mean $\varepsilon$ independent of
 $\cs$. 

We  first  study   the  process  $t\mapsto  N_{\varepsilon,t}$. 

\begin{lem}
   \label{lem:Net=Poisson}
The process $t\mapsto N_{\varepsilon,t}$ is a Poisson process
with intensity $ \displaystyle \frac{\phi_1(\gamma)}{\varepsilon
  \psi_0(\gamma)}$, where $\gamma=\psi^{-1}(1/\varepsilon)$. 
\end{lem}

\begin{proof}
We have, by the similar computations as in the  proof of Lemma \ref{lem:ccea},
\begin{align*}
\P(m_\tau= 0)
&=\inv{\varepsilon} \E\left[\int_0^\infty dt\;  \expp{-t/\varepsilon}
  \ind_{\{m_t=0\}}  \right]\\
&=\inv{\varepsilon} \int_0^\infty ds\;  \expp{-\gamma s}
\N\left[\int_0^\sigma dt\;  \expp{-t/\varepsilon} \ind_{\{m_t= 0\}}
\right] \\
&=\inv{\varepsilon\gamma } 
\N\left[\int_0^\sigma dt\;  \expp{-t/\varepsilon} \ind_{\{m_t= 0\}}
\right].
\end{align*}
By time reversibility and using optional projection and \reff{eq:N_s},
we have
\begin{align*}
\N\left[\int_0^\sigma dt\;  \expp{-t/\varepsilon} \ind_{\{m_t= 0\}}
\right]
&=\N\left[\int_0^\sigma dt\;  \expp{-(\sigma-t)/\varepsilon} \ind_{\{m_t= 0\}}
\right]   \\
&=\N\left[\int_0^\sigma dt\;  \expp{-\gamma \langle \rho_t, 1  \rangle }
  \ind_{\{m_t= 0\}} 
\right].
\end{align*}
The proof of  Lemma \ref{lem:A_s=0}, see
\reff{eq:calcul-v} and \reff{eq:calc-v}, gives that $\displaystyle
\P(m_\tau= 0)= \inv{\varepsilon\psi_0(\gamma) }$. Since
$\varepsilon^{-1}= \psi(\gamma)= \psi_0(\gamma) - \phi_1(\gamma)$, we get
$\displaystyle  \inv{\varepsilon} \P(m_\tau\neq 0)=
\frac{\phi_1(\gamma)}{\varepsilon 
  \psi_0(\gamma)}$. 
\end{proof}
We then get the following Corollary.
\begin{cor}
\label{cor:cvNetg}
There exists a sub-sequence $(\varepsilon_j, j\in \N)$ decreasing to $0$,
s.t. $\P$-a.s. for  any $t_0\geq 0$ and any  continuous  function
$h$ defined  on $\R_+\times \cm_f(\R_+)$ such that   $h(u, \mu)=0$  for $u\geq  t_0$, we  have, with
$\gamma_j=\psi^{-1}(1/\varepsilon^j)$,
\[
\lim_{j\rightarrow \infty } 
\phi_1(\gamma_j)^{-1} \sum_{k=1}^{\infty }
h(A_{S^{\varepsilon_j}_k}, \rho_{S_k^{\varepsilon_j}}^-) = \int_0^\infty 
h(u,\tilde\rho_u)\; du.
\]
\end{cor}

\begin{proof}
 Notice that
  as a direct consequence  of \reff{eq:psi/l} and \reff{eq:lim-phi1}, we
  get
\[
\lim_{\varepsilon \rightarrow 0}  \varepsilon \psi_0(\gamma)=1.
\]
Recall that $(A^\varepsilon_{S^\varepsilon_k}, k\geq 1)$ are the jumping
time of the Poisson  process $t\mapsto N_{\varepsilon,t}$ with parameter
$\phi_1(\gamma)/  \varepsilon   \psi_0(\gamma)$.   Standard  results  on
Poisson process implies the  vague convergence in distribution (see also
Lemma    XI.11.1    in    \cite{dvj:itpp2})   of    $\phi_1(\gamma)^{-1}
\sum_{k=1}^\infty   \delta_{A^{\varepsilon}_{S^{\varepsilon}_k}}   (dr)$
towards the  Lebesgue measure  on $\R_+$ as  $\varepsilon$ goes  down to
$0$.   Since  the  limit  is  deterministic, the  convergence  holds  in
probability  and a.s. along  a decreasing  sub-sequence $(\varepsilon_j,
j\in \N)$. 
In particular, if $g$ is a continuous function on $\R_+$ with compact
support (hence bounded), we have that a.s. 
\[
\lim_{j \rightarrow \infty }
\phi_1(\gamma_j)^{-1} \sum_{k=1}^\infty 
g(A^{\varepsilon_j}_{S^{\varepsilon_j}_k}) = \int_0^\infty 
g(u)\; du.
\]

Notice  that $A^\varepsilon_s\geq  A_s$ and  that  a.s. $A^\varepsilon_s
\rightarrow A_s$  as $\varepsilon$ goes  down to $0$. This  implies that
a.s.  $(A^\varepsilon_s, s\geq  0)$  converges uniformly  on compacts  to
$(A_s,s\geq 0)$.  Therefore, if  $g$ is continuous with compact support,
we have a.s.
$$\lim_{j\to+\infty}\phi_1(\gamma_j)^{-1}\sum_{k=1}^\infty\bigl|g(A_{S_k^{\varepsilon_j}}^{\varepsilon_j})-g(A_{S_k^{\varepsilon_j}})\bigr|=0.$$

So we have that
\begin{equation}\label{eq:conv}
\lim_{j \rightarrow \infty }
\phi_1(\gamma_j)^{-1} \sum_{k=1}^\infty 
g(A_{S^{\varepsilon_j}_k}) = \int_0^\infty 
g(u)\; du
\end{equation}
and this convergence also holds for a c\`ad-l\`ag function $g$ with
compact support as the Lebesgue measure does not charge the point of
discontinuity of $g$.

Let $h$ be a continuous function defined  on $\R_+\times \cm_f(\R_+)$
such that   $h(u, \mu)=0$  for $u\geq  t_0$. First let us remark that $\rho^-_{S_k^\varepsilon}    =
\rho_{T_k^\varepsilon}$  and  that  $m_{T^\varepsilon_k}=0$.  Using  the
strong Markov property at time $T_k^\varepsilon$ and the second part of
Corollary \ref{cor:C_0=0}, we deduce that $\P$-a.s. for all $k\in
\N^*$, 
\begin{equation}
    \label{eq:cat=t}
C_{A_{T_k^\varepsilon}}=T_k^\varepsilon. 
\end{equation}
Therefore, as
$A_{S_k^\varepsilon}= A_{T_k^\varepsilon}$, 
we have
$\P$-a.s.
\[
\tilde \rho_{A_{S_k^\varepsilon}}=
\tilde \rho_{A_{T_k^\varepsilon}}=
\rho_{T_k^\varepsilon}= \rho_{S_k^\varepsilon}^-.
\]

This gives
$$\phi_1(\gamma_j)^{-1} \sum_{k=1}^{\infty }
h(A_{S^{\varepsilon_j}_k}, \rho_{S_k^{\varepsilon_j}}^-)=\phi_1(\gamma_j)^{-1} \sum_{k=1}^{\infty }
h(A_{S^{\varepsilon_j}_k},\tilde \rho_{A_{S_k^{\varepsilon_j}}})$$
and applying the convergence \reff{eq:conv} to the c\`ad-l\`ag function
$$g(u)=h(u,\tilde\rho_u)$$
gives the result of the lemma.
\end{proof}

We now study $K_\varepsilon$ given by \reff{eq:Ke}. We keep the same
notation as in Lemma \ref{lem:ccea}. 
\begin{lem}
   \label{lem:CV-Ke}
There exists
a deterministic function $\mathcal R$ s.t. $\displaystyle
\lim_{\varepsilon\rightarrow 0} \mathcal R(\varepsilon) =0$ and  for all
$\varepsilon>0$ and $\mu\in\cm_f(\R_+)$, we have: 
\[
\sup_{r\geq 0}\left|\phi_1(\gamma)
\log(K_\varepsilon(r,\mu)) - \alpha_1 v(r,\mu) - \int_{(0, \infty )}
  \pi_1(d\ell_1)\;   \left(1- w(\ell_1,r,\mu) \right) \right| \leq \mathcal R(\varepsilon).
\]
\end{lem}

\begin{proof}
   We have  
\begin{align*}
   K_\varepsilon(r, \mu)
&= \frac{\psi(\gamma)}{ \psi(\gamma)- \psi(v(r,\mu))
}\frac{\gamma - v(r,\mu) }{\gamma} \inv{\phi_1(\gamma)}  \\
&\hspace{2cm}
 \left( \alpha_1 \gamma +   
 \gamma \int_0^1 \!\!du  \!  \int_{(0, \infty )}\! \!\!\! \!\!\ell_1
 \pi_1(d\ell_1)\; 
  w(u\ell_1,r,\mu) \expp{-\gamma(1-u)\ell_1} \right)\\
&= \frac{\psi(\gamma)}{ \psi(\gamma)- \psi(v(r,\mu))
}\frac{\gamma - v(r,\mu) }{\gamma} \inv{\phi_1(\gamma)}  \\
&\hspace{2cm}
 \left( \alpha_1 \gamma +   
 \int_{(0, \infty )}  \pi_1(d\ell_1)\; \int_0^{\gamma \ell_1} \expp{-s}
 ds\; 
  w(\ell_1- \frac{s}{\gamma},r,\mu)  \right)\\
&= \frac{\psi(\gamma)}{ \psi(\gamma)- \psi(v(r,\mu))
}\frac{\gamma - v(r,\mu) }{\gamma} \inv{\phi_1(\gamma)} \\
&\hspace{2cm}
 \left( \phi_1( \gamma) -
 \int_{(0, \infty )}  \pi_1(d\ell_1)\; \int_0^{\gamma \ell_1} \expp{-s}
 ds\; 
  \left(1- w(\ell_1- \frac{s}{\gamma},r,\mu) \right) \right). 
\end{align*}
In particular, we have $\phi_1(\gamma)
\log(K_\varepsilon(r,\mu))=-A_1+A_2+A_3$, where 
\begin{align*}
   A_1(r)&= \phi_1(\gamma)\log\Big(1- \psi(v(r,\mu))/\psi(\gamma)\Big),\\
   A_2(r)&= \phi_1(\gamma)\log(1- v(r,\mu)/\gamma),\\
   A_3(r)&= \phi_1(\gamma)\log\left(1-  \int_{(0, \infty )}
     \pi_1(d\ell_1)\; \int_0^{\gamma \ell_1} \expp{-s}  ds\;  
  \left(1- w(\ell_1- \frac{s}{\gamma},r,\mu) \right) /\phi_1(\gamma)\right). 
\end{align*}
Thanks to $(h_3)$, there exists a finite constant $a>0$ 
s.t. $\P$-a.s. $v(r,\mu)<a$ for all $r\geq 0$. We deduce there exists
$\varepsilon_0>0$ and a finite constant $c>0$ s.t.
$\P$-a.s  for all $\varepsilon\in (0, \varepsilon_0]$, 
\begin{equation}
   \label{eq:majo-A1}
\sup_{r\geq 0} | A_1(r) |\leq c
\frac{\phi_1(\gamma)}{\psi(\gamma)}\quad\text{and}\quad  
\sup_{r\geq 0}|A_2(r) - \alpha_1 v(r,\mu)|\leq \frac{c}{\gamma}+
c|\frac{\phi_1(\gamma)}{\gamma} - \alpha_1|. 
\end{equation}
We have 
\begin{multline*}
   \int_{(0, \infty )}  \pi_1(d\ell_1)\;
\int_0^{\gamma \ell_1} \expp{-s}  ds\;  
  \left(1- w(\ell_1- \frac{s}{\gamma},r,\mu) \right) - \int_{(0, \infty )}
  \pi_1(d\ell_1)\;   \left(1- w(\ell_1,r,\mu) \right) \\
= \int_{(0, \infty )}  \pi_1(d\ell_1)\;
 \expp{-\gamma \ell_1}   (w(\ell_1,r,\mu) - 1)\\
+  \int_{(0, \infty )}  \pi_1(d\ell_1)\;
\int_0^\infty \expp{-s}  ds\;  
  \left(w(\ell_1,r,\mu) - w(\ell_1- \frac{s}{\gamma},r,\mu) \right) \ind_
  {\{s\leq \gamma 
    \ell_1\}} .
\end{multline*}
It is then easy to get, using $(h_3)$ and \reff{eq:c-eta},  that $\P$-a.s 
\[
\phi_2(\gamma)=\sup_{r\geq 0}\left| \int_{(0, \infty )} \!\! \!\! \!\!
  \pi_1(d\ell_1)\; 
\int_0^{\gamma \ell_1} \expp{-s}  ds\;  
  \left(1- w(\ell_1- \frac{s}{\gamma},r,\mu) \right) - \int_{(0, \infty )}
   \!\! \!\! \!\! \pi_1(d\ell_1)\;   \left(1- w(\ell_1,r,\mu) \right) \right|\\
\]
converges to $0$ as $\gamma$ goes to infinity.

Recall that we assumed  that $\displaystyle \lim_{\gamma\rightarrow\infty }
\phi_1(\gamma)=+\infty $. Thus, there exist
$\varepsilon_0>0$ and a finite constant $c>0$ s.t.
$\P$-a.s  for all $\varepsilon\in (0, \varepsilon_0]$, 
\begin{equation}
   \label{eq:majo-A3}
\sup_{r\geq 0}\left|A_3(r) - \int_{(0, \infty )}
  \pi_1(d\ell_1)\;   \left(1- w(\ell_1,r,\mu) \right) \right|
\leq   \frac{c}{\phi_1(\gamma)} + \phi_2(\gamma).
\end{equation}
Using \reff{eq:majo-A1} and \reff{eq:majo-A3}, we get that there exists
a deterministic function $\mathcal R$ s.t. $\P$-a.s 
\[
\sup_{r\geq 0}\left|\phi_1(\gamma)
\log(K_\varepsilon(r,\mu)) - \alpha_1 v(r,\mu) - \int_{(0, \infty )}
  \pi_1(d\ell_1)\;   \left(1- w(\ell_1,r,\mu) \right) \right| \leq \mathcal R(\varepsilon),
\]
where $\displaystyle  \lim_{\varepsilon\rightarrow 0} \mathcal R(\varepsilon) =0$,
thanks to 
\reff{eq:psi/l} and \reff{eq:lim-phi1}.

\end{proof}

The previous results allow us to compute the following limit. 
We keep the same
notation as in Lemma \ref{lem:ccea}. 
\begin{lem}
   \label{lem:comp_lim}
Let $\varphi$ satisfying condition
   ($h_1$)--($h_3$). There exists a sub-sequence $(\varepsilon_j, j\in
   \N)$ decreasing to $0$, s.t. $\P$-a.s.
\[
\lim_{j\rightarrow\infty } \prod_{k=1}^\infty 
K_{\varepsilon_j}(A_{S^{\varepsilon_j}_k},\rho_{S_k^\varepsilon}^-)= \exp -
\int_0^\infty  
du\; \left(\alpha_1 v(u) + \int_{(0,\infty )}\pi_1(d\ell)\; (1-
  w(\ell,u,\mu))\right).
\]
\end{lem}

\begin{proof}
   Notice that thanks to ($h_1$), the functions $v$ and
   $(u,\mu)\mapsto w(\ell,u,\mu)$
   are continuous and that for $r\geq t$, $v(r,\mu)=0$ and
   $w(\ell,r,\mu)=1$. The result is then a direct consequence of Corollary 
   \ref{cor:cvNetg} and Lemma \ref{lem:CV-Ke}. 
\end{proof}

\subsection{Proof of Theorem \ref{th:SMP}}
\label{sec:proof_th_MS}
Now we can prove the special Markov property in the case
$\lim_{\gamma\rightarrow\infty }\phi_1(\gamma)=+\infty $.

Let $Z\in  \tilde \cf_\infty $  non-negative such that  $\E[Z]<\infty $.
Let $\varphi$  satisfying hypothesis  of Theorem \ref{th:SMP},
($h_1$)--($h_3$). We have, using notation of the previous sections
\begin{align*}
 \E\left[Z \exp\left(-\sum_{i\in
    I}\varphi(A_{\alpha_i},\rho_{\alpha_i-},\cs^i)\right)\right]
&=  \lim_{j\rightarrow\infty } \E\left[ Z \exp\left(
    -\sum_{k=1}^\infty \varphi^*\left
      (A_{S_k^{\varepsilon_j} } ,\rho_{S_k^{\varepsilon_j}}^-,\cs^{k,\varepsilon_j }
    \right) \right) 
\right] \\ 
&=  \lim_{j\rightarrow\infty } \E\left[ Z \prod_{k=1}^\infty 
K_{\varepsilon_j}(A_{S^{\varepsilon_j}_k},\rho_{S_k^{\varepsilon_j}}^-)\right] \\ 
&= \E\left[ Z \expp{ - \int_0^\infty 
du\; \left(\alpha_1 v(u,\tilde\rho_u) + \int_{(0,\infty )}\pi_1(d\ell)\; (1-
  w(\ell,u,\tilde\rho_u))\right)}\right] , 
\end{align*}
where we  used Lemma \ref{lem:approx} and dominated  convergence for the
first  equality,  Lemma \ref{lem:ccea}  for  the  second equality,  Lemma
\ref{lem:comp_lim}  and dominated
convergence for the last equality. By 
monotone class Theorem and  monotonicity, we can remove hypothesis
($h_1$)-- ($h_3$). To ends the proof of the first part, notice that 
$\int_0^t
du\; \left(\alpha_1 v(u) + \int_{(0,\infty )}\pi_1(d\ell)\; (1-
  w(\ell,u))\right)$ is $\tilde \cf_\infty$-measurable and so this is
$\P$-a.e. equal to the conditional expectation (i.e. the left hand side
term of \reff{eq:MS}). 

\section{Law of the pruned exploration process}
\label{sec:law_pruned}
Let $\rho^{(0)}$ be the 
exploration  process of a L\'evy process with Laplace exponent $\psi_0$. 
The  aim of this section is to prove  Theorem \ref{thm:law_pruned}. 


\subsection{A martingale problem for $\tilde \rho$}
Let $\tilde\sigma=\inf\{t>0,\tilde\rho_t=0\}$.
In this  section, we  shall compute  the law of  the total  mass process
$(\langle \tilde\rho_{t\wedge \tilde  \sigma} ,1\rangle,\ t\ge 0)$ under
$\P_\mu=\P_{\mu,0}$, using martingale problem characterization.  We will
first show how a martingale problem  for $\rho$ can be translated into a
martingale problem  for $\tilde \rho$, see also \cite{ad:fpigep}.
Unfortunately, we  were not able 
to  use standard  techniques  of  random time  change,  as developed  in
Chapter 6 of  \cite{ek:mp} and used for Poisson  snake in \cite{as:psf},
mainly     because     $\displaystyle     t^{-1}    \left(\E_{\mu}     [
  f(\rho_t)\ind_{\{m_t=0\}}] -  f(\mu)\right)$ may not have  a limit as
$t$ goes down to 0, even for exponential functionals.

Let  $F,  K\in  \cb(\cm_f(\R_+))$ be bounded. We suppose that  $\displaystyle
\N\left[\int_0^\sigma \val{K(\rho_s)}\; 
  ds  \right]<\infty  $, that  for  any  $\mu\in
\cm_f(\R_+)$,  $\displaystyle  \E_\mu^*\left[\int_0^\sigma \val{K(\rho_s)}\;
  ds  \right]<\infty  $  and  that $M_t=F(\rho_{t\wedge  \sigma})
-\int_0^{t\wedge   \sigma}  K(\rho_s)\; ds$,   for  $t\geq   0$,   defines  an
$\cf$-martingale. In other words, if $F$ belongs to the domain of the
infinitesimal generator $\cl$ of $\rho$, we have $K=\cl F$. We will
see in the proof of Corollary \ref{cor:law_total_mass} that these assumptions on $F$ and $K$ are in
particular fulfilled for
$$F(\nu)=\expp{-c\langle\nu,1\rangle}\quad\mbox\quad K(\nu)=\psi(c)F(\nu).$$

    Notice     that we have
$$|M_t|\le \|F\|_\infty+\int_0^\sigma \bigl|K(\rho_s)\bigr|ds$$
and thus   $\E_\mu^*\left[\sup_{t\geq    0}
  \val{M_t}\right]<\infty $. Consequently, we can define for $t\geq 0$,
\[
N_t=\E^*_\mu[M_{C_t}|\tilde \cf_t].
\]
\begin{prop}
\label{prop:mart-tr}
 The process  $N=(N_t, t\geq 0)$ is an $\tilde \cf$-martingale.
And we have the representation formula for $N_t$: 
\begin{equation}
   \label{eq:Mart-tr}
N_t=F(\tilde \rho_{t\wedge \tilde \sigma})  - 
\int_0^{t\wedge  \tilde \sigma} du\; \tilde  K(\tilde \rho_u) ,
\end{equation}
with
\begin{equation}
   \label{eq:barK}
\tilde K(\nu)=  K(\nu) +
\alpha_1 \N\left[\int_0^\sigma {K([\nu, \rho_s])} \; ds\right]+
\int_{(0,\infty )} \pi_1(d\ell)
\;\E_\ell ^*
\left[\int_0^\sigma {K([\nu, \rho_s])} \; ds\right].
\end{equation}
\end{prop}

\begin{proof}
  Notice that $N=(N_t, t\geq 0)$ is an $\tilde \cf$-martingale. Indeed, we have
for $t,s\geq 0$, 
\begin{align*}
   \E_{\mu}[N_{t+s}|\tilde \cf_t]
&=\E_{\mu}[\E_{\mu}[M_{C_{t+s}}|\tilde \cf_{t+s}]|\tilde \cf_t]\\
&=\E_{\mu}[M_{C_{t+s}}|\tilde \cf_t]\\
&=\E_{\mu}[\E_{\mu}[M_{C_{t+s}}|\cf_{C_t}]|\tilde \cf_t]\\
&=\E_{\mu}[M_{C_{t}}|\tilde \cf_t],
\end{align*}
where we used the optional stopping time Theorem for the last equality.
To compute $\E_{\mu}[M_{C_t}|\tilde \cf_t]$, we set
$N'_t =M_{C_t}+
M'_{C_t} $, where for $u\geq 0$, 
\[
M'_{u}=\int_0^{u\wedge \sigma} K(\rho_s)\ind_{\{m_s\neq 0\}} \; ds.
\]

Recall that $C_0=0$ $\P_{\mu}$-a.s. by Corollary \ref{cor:C_0=0}.
In particular, we get
\begin{align*}
   N'_t= M_{C_t} + M'_{C_t}
&=F(\rho_{C_t\wedge \sigma})  - 
\int_0^{C_t\wedge \sigma} K(\rho_s) \ind_{\{m_s= 0\}}  \; ds \\
&=F(\tilde \rho_{t\wedge \tilde \sigma})  - 
\int_0^{C_t\wedge \sigma} K(\rho_s)  \; dA_s \\
&=F(\tilde \rho_{t\wedge \tilde \sigma})  - 
\int_0^{t\wedge  \tilde \sigma} K(\tilde \rho_u)  \; du ,
\end{align*}
where  we used  the  time change  $u=A_s$  for the  last  equality.   In
particular, as  $\tilde \sigma $ is  an $\tilde \cf$-stopping time,  we get that
the process $(N'_t, t\geq 0)$  is $\tilde\cf$-adapted.  Since  $N_t=N'_t
-  \E_{\mu}[M'_{C_t}|\tilde \cf_t]  $, we  are left  with  the  computation  of
$\E_{\mu}[M'_{C_t}|\tilde \cf_t] $.

We keep  the notations of Section  \ref{sec:Markov_special}. We consider
$(\rho^i,  m^i)$, $i\in  I$  the excursions  of  the process  $(\rho,m)$
outside  $\{s,m_s=0\}$ before  $\sigma$ and  let  $(\alpha_i, \beta_i)$,
$i\in I$ be the corresponding  interval excursions. In particular we can
write
\[
\int_0^{C_t\wedge \sigma} \val{K(\rho_s)} \ind_{\{m_s\neq 0\}}\;
ds  =\sum_{i\in    I}     \Phi( A_{\alpha_i},
\rho_{\alpha_i-},\rho^i) ,
\]
with
\[
\Phi(u,\mu,\rho)=\ind_{\{u<t\}} \int_0^{\sigma(\rho)}  \val{K([\mu, \rho_s])} \; ds,
\]
where $\sigma(\rho)=\inf\{v>0; \rho_v=0\}$. 
We deduce from the second part of Theorem \ref{th:SMP},  that $\P_\mu$-a.s. 
\begin{equation}
    \label{eq:K-hat}
\E_{\mu}\left[\int_0^{C_t\wedge \sigma} \val{K(\rho_s)} \ind_{\{m_s\neq 0\}}\;
ds|\tilde \cf_\infty   \right]  =   \int_0^{\tilde   \sigma} 
\ind_{\{u<t\}} \hat K(\tilde \rho_u) \; du ,
\end{equation}
with, $\hat K$ defined for $\nu\in \cm_f(\R_+)$ by  
\[
\hat K(\nu) 
= \alpha_1 \N\left[\int_0^\sigma \val{K([\nu, \rho_s])} \; ds\right]+ 
\int_{(0,\infty )} \pi_1(d\ell)
\;\E_\ell ^*
\left[\int_0^\sigma \val{K([\nu, \rho_s])} \; ds\right].
\]
Since  $\E_\mu\left[\int_0^{C_t \wedge \sigma }
  \val{K(\rho_s)}\ind_{\{m_s \neq 0\}} \;
  ds  \right]\leq \E_\mu\left[\int_0^{\sigma }
  \val{K(\rho_s)} \;
  ds  \right]
<\infty  $, we deduce from \reff{eq:K-hat}  that $\P_{\mu}$-a.s.  $du$-a.e.   $\ind_{\{u<\tilde
  \sigma\}}\hat K(\tilde \rho_u) $ is finite.

We  define $\bar   K\in \cb(\cm_f(\R_+))$  for $\nu\in  \cm_f(\R_+)$ by
\begin{equation}
   \label{eq:def-tK}
\bar   K(\nu) 
= \alpha_1 \N\left[\int_0^\sigma {K([\nu, \rho_s])} \; ds\right]+ 
\int_{(0,\infty )} \pi_1(d\ell)
\;\E_\ell ^*
\left[\int_0^\sigma {K([\nu, \rho_s])} \; ds\right]
\end{equation}
if   $\hat  K(\nu)<\infty   $,  or   by  $\bar  K(\nu)=0$   if  $\hat
K(\nu)=+\infty  $.  In  particular,  we have  $|\bar  K(\nu)|\leq  \hat
K(\nu)$  and $\P_{\mu}$-a.s.   $\int_0^{\tilde \sigma}  |\bar  K(\tilde
\rho_u)| \;  du$ is finite.  Using the special Markov property once again
(see \reff{eq:K-hat}), we get that $\P_\mu$-a.s.,
\[
\E_{\mu}\left[M'_{C_t} |\tilde \cf_\infty   \right]  
= \E_{\mu}\left[\int_0^{C_t\wedge \sigma} {K(\rho_s)} \ind_{\{m_s\neq 0\}}\;
ds|\tilde \cf_\infty   \right]  =   \int_0^{t\wedge\tilde   \sigma} 
 \bar K(\tilde \rho_u) \; du. 
\]
Finally, as $N_t
=N'_t- \E_{\mu}\left[M'_{C_t} |\tilde \cf_\infty   \right] $,  this
gives \reff{eq:Mart-tr}. 
\end{proof}

\begin{cor}\label{cor:law_total_mass}
  Let $\mu\in \cm_f(\R_+)$. The law  of the total mass process $(\langle
  \tilde\rho_t,1\rangle,\ t\ge  0)$ under  $\P_{\mu,0}^*$ is the  law of
  the total mass process of $\rho^{(0)}$ under $\P_\mu^*$.
\end{cor}

\begin{proof}
  Let $X=(X_t, t\geq 0)$ be under $\rP^*_x$ a L\'evy process with Laplace
  transform $\psi$  started at  $x>0$ and stopped  when it  reached $0$.
  Under  $\P_\mu$,  the   total  mass  process  $(\langle  \rho_{t\wedge
    \sigma},  1\rangle,   t  \geq  0)$  is  distributed   as  $X$  under
  $\rP^*_{\langle \mu,  1\rangle} $.  Let  $c> 0$.  From  L\'evy processes
  theory, we  know that  the process $\expp{-c  X_t} -  \psi(c) \int_0^t
  \expp{-cX_s}\; ds $, for $t\geq 0$ is a martingale. We deduce from the
  stopping time  Theorem that $M=(M_t, t\geq 0)$  is an $\cf$-martingale
  under  $\P_\mu$,  where $M_t=F(\rho_{t\wedge  \sigma})-\int_0^{t\wedge
    \sigma} K(\rho_s) \; ds $,  with $F, K\in \cb( \cm_f(\R_+))$ defined
  by  $F(\nu)=\expp{   -  c  \langle   \nu,  1\rangle}$  for   $\nu  \in
  \cm_f(\R_+)$ and $K=\psi(c)F$. Notice  $K\geq 0$. We have by dominated
  convergence and monotone convergence.
\[
\expp{-   c\langle   \mu,   1\rangle}   =\lim_{t\rightarrow   \infty   }
\E_\mu[M_t]=\E_\mu[\expp{  - c  \langle \rho_\sigma  ,  1\rangle}] -
\psi(c)  \E_\mu  \left[\int_0^\sigma  \expp{  -  c  \langle  \rho_s  ,
    1\rangle} \; ds \right] .
\]
This implies that, for any $\mu\in
\cm_f(R_+)$,  $\displaystyle  \E_\mu\left[\int_0^\sigma \val{K(\rho_s)}\;
  ds  \right]$  is finite. Using the Poisson representation, see
Proposition  \ref{prop:poisson_representation2}, it is easy to get that 
\begin{equation}
   \label{eq:Necrho}
\N\left[\int_0^\sigma dt\;  \expp{  - c  \langle \rho_t  ,
    1\rangle}\right]=\frac{c}{\psi(c)}.
\end{equation}
In particular, it is also finite.

{F}rom Proposition \ref{prop:mart-tr}, we get that $N=(N_t, t\geq 0)$
is under $\P_\mu$ an $\tilde \cf$-martingale, where:  for $t\geq 0$,  
\[
N_t=\expp{-c   \langle \tilde\rho_{t\wedge \tilde \sigma}  ,  1\rangle} -
 \int_0^{t\wedge \tilde \sigma} 
\tilde K(\tilde \rho_u)\; du
\]
and $\tilde K$ given by \reff{eq:barK}. We can compute $\tilde K$:
\begin{align*}
   \tilde  K(\nu)
&= \psi(c) \expp{- c\langle \nu,1 \rangle} \left(1+
\alpha_1  \N\left[\int_0^\sigma
  \expp{- c\langle \rho_s,1 \rangle} \; ds\right]+
\int_{(0,\infty )} \pi_1(d\ell)
\;\E_\ell ^*
\left[\int_0^\sigma  \expp{- c\langle \rho_s,1 \rangle} \;
  ds\right]\right)\\
&= \psi(c) \expp{- c\langle \nu,1 \rangle} \left(1+
\alpha_1 \frac{c}{\psi(c)} +
\int_{(0,\infty )} \pi_1(d\ell) \int_0^\ell dr\; \expp{-cr }
\N
\left[\int_0^\sigma  \expp{- c\langle \rho_s,1 \rangle} \;
  ds\right]\right)\\
&= \expp{- c\langle \nu,1 \rangle} \left(\psi(c) + \alpha_1 c +
  \int_{(0,\infty )} \pi_1(d\ell) (1- \expp{-c\ell })\right)\\
&=\psi_0(c)\expp{- c\langle \nu,1 \rangle} ,
\end{align*}
where we  used \reff{eq:Necrho} and the excursion  decomposition for the
second equality, and $\psi_0=\psi+\phi_1$ for the last one.

Thus, the process $(N_t, t\geq 0)$ with for $t\geq 0$ 
 \[
N_t=\expp{-c   \langle \tilde\rho_{t\wedge \tilde \sigma}  ,  1\rangle} -
\psi_0(c)  \int_0^{t\wedge \tilde \sigma} 
\expp{-c   \langle \tilde\rho_{u}  ,  1\rangle}\; du
\]
is under $\P_\mu$ an $\tilde \cf$-martingale. 

Notice  that $\tilde  \sigma=\inf  \{ s\geq  0;  \langle \tilde\rho_s  ,
1\rangle=0\}$.  Let $X^{(0)}=(X^{(0)}_t ,  t\geq 0)$  be under
$\rP^*_x$  a  L\'evy  process  with Laplace  transform  $\psi_0$
started at $x>0$ and stopped when it reached $0$.
The two non-negative  c\`ad-l\`ag processes $( \langle
  \tilde\rho_{t\wedge \tilde \sigma}   ,  1\rangle, t\geq 0)$ and 
$X^{(0)}$ solves the martingale problem: for any $c\geq 0$, the
process defined for $t\geq 0$ by 
\[
\expp{- c Y_{t\wedge \sigma'}} - \psi_0(c) \int_0^{t\wedge \sigma'} \expp{- c Y_s}\; ds,
\]
where $\sigma'=\inf  \{ s\geq  0;  Y_s\leq 0\}$, is a martingale. {F}rom
Corollary 4.4.4 in  \cite{ek:mp}, we deduce that those two processes have
the same distribution. To finish the proof, notice that the total mass
process of $\rho^{(0)}$ under $\P_\mu^*$ is distributed as
$X^{(0)}$ under $\rP^*_{\langle \mu,  1\rangle} $.
\end{proof}

\subsection{Identification of the law of $\tilde \rho$}

To begin with, let us mention some useful properties of the process $\tilde\rho$.
\begin{lem}
   We have the following properties for the process $\tilde \rho$. 
\begin{enumerate}
   \item[(i)] $\tilde \rho$ is a c\`ad-l\`ag Markov process.
   \item[(ii)] The sojourn time at $0$ of $\tilde \rho$ is  $0$.
   \item[(iii)] 0 is recurrent for $\tilde \rho$.
\end{enumerate}
\end{lem}

\begin{proof}
   
(i) This is a direct consequence of the strong Markov property of the
    process $(\rho, m)$.

(ii) We have for $r>0$, with the change of variable $t=A_s$, a.s. 
\[
\int_0^r \ind_{\{\tilde \rho _t=0\}} \; dt
=\int_0^r \ind_{\{ \rho _{C_t}=0\}} \; dt
=\int_0^{C_r} \ind_{\{\rho _s=0\}} \; dA_s
=\int_0^{C_r} \ind_{\{\rho _s=0\}} \; ds=0,
\] 
as the sojourn time of $\rho$ at $0$ is $0$ a.s.

(iii) Since $ \tilde  \sigma =A_\sigma$ and $\sigma<+\infty $ a.s., we
    deduce that $0$ is recurrent for $\tilde \rho$ a.s. 
\end{proof}

Since the processes $\tilde \rho$ and $\rho^{(0)}$ are both
Markov processes, to show that they have the same law, it is enough to
show that they have the same one-dimensional marginals. We first prove
that result under the excursion measure.

\begin{prop}\label{prop:=n}
For every $\lambda>0$ and every non-negative bounded measurable
function $f$,
$$\N\left[\int_0^{\tilde\sigma}\expp{-\lambda t-\langle \tilde\rho_t,f\rangle}dt\right]=\N\left[\int_0^{\sigma^{(0)}}\expp{-\lambda t-\langle \rho^{(0)}_t,f\rangle}dt\right].$$
\end{prop}

\begin{proof}
On one hand, we compute, using the definition of the pruned process
$\tilde \rho$,
$$
\N\left[\int_0^{\tilde\sigma}\expp{-\lambda t-\langle
    \tilde\rho_t,f\rangle}dt\right]
=\N\left[\int_0^{A_\sigma}\expp{-\lambda t-\langle \rho_{C_t},f\rangle}dt\right].
$$
We now make the change of variable $t=A_u$ to get
\begin{align*}
\N\left[\int_0^{\tilde\sigma}\expp{-\lambda t-\langle
    \tilde\rho_t,f\rangle}dt\right] & =\N\left[\int_0^\sigma
    \expp{-\lambda A_u}\expp{-\langle\rho_u,f\rangle}dA_u\right]\\
& =\N\left[\int_0^\sigma \expp{-\lambda A_u}\expp{-\langle\rho_u,f\rangle}\ind_{\{m_u=0\}}du\right].
\end{align*}
By a time reversibility argument, see Lemma \ref{lem:reversib}, we obtain
\begin{align*}
\N\left[\int_0^{\tilde\sigma}\expp{-\lambda t-\langle
    \tilde\rho_t,f\rangle}dt\right]
& =\N\left[\int_0^\sigma
    \ind_{\{m_u=0\}}\expp{-\langle\eta_u,f\rangle}\expp{-\lambda(A_\sigma-A_u)}du\right]\\
& =\N\left[\int_0^\sigma
    \ind_{\{m_u=0\}}\expp{-\langle\eta_u,f\rangle}\E_{\rho_u,0}^*\left[\expp{-\lambda
    A_\sigma}\right]du\right]\\
& =\N\left[\int_0^\sigma
    \ind_{\{m_u=0\}}\expp{-\langle\eta_u,f\rangle
      -{\psi_0^{-1}(\lambda)} \langle 
      \rho_u,1\rangle}du\right], 
\end{align*}
where we applied Lemma \ref{lem:A_s=0} (i) for the last equality. 
Now, using Proposition \ref{prop:poisson_representation2}, we have 
\[
\N\left[\int_0^{\tilde\sigma}\expp{-\lambda t-\langle
    \tilde\rho_t,f\rangle}dt\right]
=\int_0^\infty  da \; \expp{-\alpha a}  \M_a \left[
  \ind_{\{m=0\}}\expp{-\langle \nu,f\rangle 
      -{\psi_0^{-1}(\lambda)} \langle 
      \mu,1\rangle}  \right].
\]
Using usual properties of point Poisson measures, we have, with $c=\alpha_1+\int_{(0,\infty
  )}\ell\; 
\pi_1(d\ell)$, 
\[
\M_a \left[ \ind_{\{m=0\}} F(\mu,\nu) \right]
=\expp{- c a  }\M_a \left[  F(\mu^0,\nu^0) \right],
\]
where with the notations of Proposition
\ref{prop:poisson_representation2}, for any $f\in \cb_+(\R_+)$  
\begin{align*}
\langle \mu^0,f\rangle  & =\int\mathcal{N}_0(dx\, d\ell\, 
du)\, \ind_{[0,a]}(x) u\ell f(x) + \beta \int_0^a f(r)\; dr,\\
\langle \nu^0,f\rangle  & =\int\mathcal{N}_0(dx\, d\ell\, 
du)  \, \ind_{[0,a]}(x) (1-u)\ell f(x)+ \beta \int_0^a f(r)\; dr.
\end{align*}
As $\alpha_0=\alpha+c$,  we have 
\[
\N\left[\int_0^{\tilde\sigma}\expp{-\lambda t-\langle
    \tilde\rho_t,f\rangle}dt\right]
=\int_0^\infty  da \; \expp{-\alpha_0 a}\M_a \left[ \expp{-\langle
    \nu^0,f\rangle        -{\psi_0^{-1}(\lambda)} \langle 
      \mu^0,1\rangle}  \right].
\]
Proposition  3.1.3  in  \cite{dlg:rtlpsbp}  directly  implies  that  the
left-hand   side  of the previous equality  is   equal   to   $\displaystyle   \N   \left[\int_0^
  {\sigma^{(0)}}\expp{-\langle\eta   ^{(0)}_t,f\rangle   -  {\psi_0^{-1}
      (\lambda)} \langle \rho_t^{(0)},1\rangle }dt\right]$.
On the other hand, similar 
computations as above yields that this quantity is equal to 
$\displaystyle  \N \left[\int_0^
  {\sigma^{(0)}}\expp{-\lambda t- \langle\rho   ^{(0)}_t,f\rangle }dt\right]$. This ends the proof. 
\end{proof}

Now, we prove the same result under $\P_{\mu,0}^*$, that is:
\begin{prop}\label{prop:Laplace_*}
For every $\lambda>0$,  $f\in \cb_+(\R_+)$ bounded and every finite measure $\mu$,
$$\E_{\mu,0}^*\left[\int_0^{\tilde\sigma}\expp{-\lambda t-\langle \tilde\rho_t,f\rangle}dt\right]=\E_\mu^*\left[\int_0^{\sigma^{(0)}}\expp{-\lambda t-\langle \rho^{(0)}_t,f\rangle}dt\right].$$
\end{prop}

\begin{proof}
{F}rom the Poisson representation, see Lemma \ref{lem:dlg-decomp}, and
using notations of this Lemma and of \reff{eq:def-Ai} we have 
\begin{align*}
\E_{\mu,0}^*\left[\int_0^{\tilde \sigma}\expp{-\lambda
    t-\langle\tilde\rho_t,f\rangle}dt\right] 
& =\E_{\mu,0}^*\left[\int_0^\sigma \expp{-\lambda A_u-\langle
    \rho_u,f\rangle}dA_u\right]\\
& =\E_{\mu,0}^*\left[\sum_{i\in J}\expp{-\lambda
    A_{\alpha_i}-\langle
    k_{-I_{\alpha_i}},f\rangle}\int_0^{\sigma_i}\expp{-\langle
    \rho_s^i,f_{-I_{\alpha_i}}\rangle-\lambda A^i_s}dA^i_s\right],
\end{align*}
where the function $f_r$ is defined by $f_r(x)=f(H^{(\mu)}_r+x)$ and 
$H^{(\mu)}_r=H(k_r\mu)$ is the maximal
element of the 
closed support of $k_r \mu$ (see \reff{def:H}).
We recall that $-I$ is the local time at 0 of the reflected process 
 $X-I$, and that $\tau_r=\inf\{s; -I_s>r\}$ is the  right continuous inverse of $-I$. {F}rom excursion formula, and using the time
change $-I_s=r$ (or equivalently $\tau_r=s$), we get 
\begin{align}
\nonumber
\E_{\mu,0}^*\left[\int_0^{\tilde \sigma}\expp{-\lambda
    t-\langle\tilde\rho_t,f\rangle}dt\right] 
& =\E_{\mu,0}^*\left[\int_0^{\tau_{\langle\mu,1\rangle}}d(-I_s)\expp{-\langle
    k_{-I_s}\mu,f\rangle-\lambda A_s}G(-I_s)\right]\\
& =\E_{\mu,0}^*\left[\int_0^{\langle \mu,1\rangle} dr\expp{-\langle
    k_r\mu,f\rangle-\lambda A_{\tau_r}}G(r)\right],
\label{eq:EGr1}
\end{align}
where the function $G(r)$ is given by
\[
G(r)=\N\left[\int_0^\sigma\expp{-\langle\rho_s,f_r\rangle-\lambda
    A_s}dA_s\right]=  
\N\left[\int_0^{\tilde\sigma}\expp{-\lambda t-\langle
    \tilde\rho_t,f_r\rangle}dt\right].
\]
The same kind of computation gives
\begin{equation}
   \label{eq:EGr2}
\E_\mu^*\left[\int_0^{\sigma^{(0)}}\expp{-\lambda t-\langle
    \rho_t^{(0)},f\rangle}dt\right]=\E\left[\int_0^{\langle
    \mu,1\rangle}dr\expp{-\langle
    k_r\mu,f\rangle-\lambda\tau_r^{(0)}}G^{(0)}(r)\right]
\end{equation}
where the function $G^{(0)}$ is defined by
\[
G^{(0)}(r)=\N\left[\int_0^{\sigma^{(0)}}\expp{-\lambda
    s-\langle\rho_s^{(0)},f_r\rangle}ds\right]
\]
and $\tau^{(0)}$ is the right-continuous inverse of the infimum
    process $-I^{(0)}$ of the L\'evy process with Laplace
    exponent $\psi_0$. 

Proposition \ref{prop:=n} says that the functions $G$ and
$G^{(0)}$ are equal. Moreover, as the total mass processes have
the same law (see Corollary \ref{cor:law_total_mass}), we know that
the proposition is true for $f$ constant. And, for $f$ constant, the
functions $G$ and
$G^{(0)}$ are also constant. Therefore, we have for $f$ constant
equal to $c\geq 0$,
\[
\E_{\mu,0}^*\left[\int_0^{\langle \mu,1\rangle}dr\expp{-c(\langle
    \mu,1\rangle-r)}\expp{-\lambda                     A_{\tau_r}}\right]
=\E\left[\int_0^{\langle                  \mu,1\rangle}dr\expp{-c(\langle
    \mu,1\rangle-r)}\expp{-\lambda  \tau_r^{(0)}}\right].
\]
As this is  true for any $c\geq 0$, uniqueness  of the Laplace transform
gives the equality
\[
\E_{\mu,0}^*\left[\expp{-\lambda
    A_{\tau_r}}\right]=\E\left[\expp{-\lambda\tau_r^{(0)}}\right]\qquad
    dr-\mbox{a.e.}
\]
In fact this equality holds for every $r$ by right-continuity.

Finally as $G=G^{(0)}$, we have thanks to \reff{eq:EGr1} and
\reff{eq:EGr2}, that, for every bounded non-negative measurable
function $f$,
$$\int_0^{\langle \mu,1\rangle} dr\expp{-\langle
    k_r\mu,f\rangle}\E_{\mu,0}^*\left[\expp{-\lambda A_{\tau_r}}\right]G(r)=\int_0^{\langle
    \mu,1\rangle}dr\expp{-\langle
    k_r\mu,f\rangle}\E\left[\expp{-\lambda\tau_r^{(0)}}\right]G^{(0)}(r)$$
which ends the proof.
\end{proof}

\begin{cor}
\label{cor:=*} 
  The  process $\tilde  \rho$ under  $\P^*_{\mu,0}$ is  distributed as
   $\rho^{(0)}$ under $\P^*_\mu$.
\end{cor}

\begin{proof}
  Let $f\in  \cb_+(\R_+)$ bounded. Proposition  \ref{prop:Laplace_*} can
  be re-written as
$$\int_0^{+\infty}\expp{-\lambda t}\E_{\mu,0}^*\left[\expp{-\langle
    \tilde\rho_t,f\rangle}\ind_{\{t\le \tilde\sigma\}}\right]dt =\int_0^{+\infty}\expp{-\lambda t}\E_{\mu}^*\left[\expp{-\langle
    \rho^{(0)}_t,f\rangle}\ind_{\{t\le\sigma^{(0)}\}}\right]dt.$$
By uniqueness of the Laplace transform, we deduce that, for almost every
    $t>0$,
$$\E_{\mu,0}^*\left[\expp{-\langle
    \tilde\rho_t,f\rangle}\ind_{\{t\le
    \tilde\sigma\}}\right]=\E_{\mu}^*\left[\expp{-\langle
    \rho^{(0)}_t,f\rangle}\ind_{\{t\le\sigma^{(0)}\}}\right].$$
In fact this  equality holds for every $t$  by right-continuity.  As the
Laplace functionals characterize the law  of a random measure, we deduce
that,  for  fixed  $t>0$, the  law of  $\tilde\rho_t$ under
$\P_{\mu,0}^*$  is the  same  as the  law  of $\rho^{(0)}_t$  under
$\P_\mu^*$.

The Markov property  then gives the equality in law for the  c\`ad-l\`ag processes
$\tilde\rho$ and $\rho^{(0)}$.
\end{proof}

\begin{proof}[Proof of Theorem \ref{thm:law_pruned}]
   $0$ is recurrent for the Markov c\`ad-l\`ag processes $\tilde \rho$ and
   $\rho^{(0)}$.  These two processes have no sojourn time at $0$, and
   when killed on the first hitting time of $0$, they have the
   same law, thanks to Lemma \ref{cor:=*}. {F}rom Theorem 4.2 of
   \cite{b:emp}, Section 5, we deduce that $\tilde \rho$ under
   $\P_{\mu,0}$ is distributed as $\rho^{(0)}$ under $\P_\mu$. 
\end{proof}

\section{Law of the excursion lengths}
\label{sec:prop_pruned}

Recall  $\tilde\sigma=\int_0^\sigma  \ind_{\{m_s=0\}}\; ds $
denotes the length of the excursion of the pruned exploration process. We
can compute  the joint  law of $(\tilde\sigma,\sigma)$.   This will
determine  uniquely  the  law  of $\tilde\sigma$  conditionally  on
$\sigma=r$.

\begin{prop}\label{prop:s-sq}
   For all non-negative $\gamma, \kappa$, the value  $v$ defined
   by $\displaystyle v=\N\left[1-\expp{-\psi(\gamma) \sigma -\kappa
   \tilde\sigma}\right] $ is the unique non-negative solution of
   the equation 
\[
\psi_0(v)= \kappa +\psi_0(\gamma).
\]
\end{prop}

\begin{proof}
Excursion theory implies that the special Markov property, Theorem
\ref{th:SMP}, also holds under $\N$, with the integration of $u$ over
$[0, \tilde \sigma=A_\sigma]$ instead of $[0,\infty )$. Taking 
   $\phi(\cs)=\psi(\gamma)\sigma$, we have
\begin{align*}
v
= \N\left[1-\expp{ -\kappa \tilde\sigma-\psi(\gamma)\sigma}\right]
&= \N\left[1-\expp{ -(\kappa+\psi(\gamma) )
    \tilde\sigma-\psi(\gamma) 
    \int_0^\sigma \ind_{\{m_s\neq  0\}}\; ds}\right]\\
&= \N\left[1-\expp{ -(\kappa+\psi(\gamma))
    \tilde \sigma- \tilde \sigma \left(\alpha_1 \N[1-\expp{-\psi(\gamma)
      \sigma}] + \int_{(0,+\infty )} \pi_1(d\ell)
    (1- \E_\ell^* [\expp{( -\psi(\gamma) \sigma)}]  \right)}\right].
\end{align*}
Notice that  $\sigma$ under $\P^*_\ell$ is distributed  as $\tau_\ell$, the
first  time  for which    the  infimum  of  $X$,  started   at  $0$,  reaches
$-\ell$. Since  $\tau_\ell$ is distributed as a
subordinator with Laplace exponent $\psi^{-1}$ at time
$\ell$,  we have 
\[
\E_\ell^* [\expp{ -\psi(\gamma) \sigma}] 
=\E\left[ \expp{-\psi(\gamma) \tau_\ell} \right] = \expp{-\ell \gamma}.
\]
Thanks to \reff{eq:N_s}, we get $\N[1-\expp{-\psi(\gamma)
      \sigma}]= \gamma$. We deduce that 
\begin{align*}
v
&= \N\left[1-\expp{ -(\kappa+\psi(\gamma))
    \tilde \sigma- \tilde \sigma \left(\alpha_1 \gamma + \int_{(0,+\infty )} \pi_1(d\ell)
    (1- \expp{ -\gamma\ell})  \right)}\right]\\
&= \N\left[1-\expp{ -(\kappa+\psi_0(\gamma))
    \tilde \sigma}\right]\\
&= {\psi_0}^{-1} (\kappa+\psi_0
    (\gamma)).
\end{align*}
Since $\psi_0$  is increasing and continuous, we get the result. 
\end{proof}

\section{Appendix}\label{sec:appendix}

We  shall  present  in  a  first  subsection, how  one  can  extend  the
construction of  the L\'evy snake from \cite{dlg:rtlpsbp}  to a weighted
L\'evy  snake, when the  height process  may not  be continuous  and the
lifetime process is  given by the total mass  of the exploration process
(instead    of   the    height   of    the   exploration    process   in
\cite{dlg:rtlpsbp}). Then,  using this construction, we can  define in a
second subsection a general L\'evy  snake when the height process is not
continuous.

\subsection{Weighted L\'evy snake}
Let $D$ be a distance on $\cm_f(\R_+)$ which defines
the topology of weak convergence. Let us recall that $(\cm_f(\R_+), D)$ 
is a Polish space, see \cite{d:mvmp}, section 3.1. 

Let  $E$ be  a  Polish space,  whose  topology is  defined  by a  metric
$\delta$, and $\partial$ be a  cemetery point added to $E$.  Let $\cw_x$
be the  space of all $E$-valued  weighted killed paths  started at $x\in
E$. An element  $\bar w= (\mu, w)$ of $\cw_x$ is  a mass measure $\mu\in
\cm_f(\R_+)$  and  a  c\`ad-l\`ag  mapping  $w:[0,  \langle  \mu,1  \rangle)
\rightarrow  E$ s.t.  $w(0)=x$.  By  convention the  point  $x$ is  also
considered as a  weighted killed path with mass  measure $\mu=0$. We set
$\cw=\bigcup _{x\in E}\cw_x$ and equip $\cw$ with the distance
\begin{multline}\label{eq:dist_d}
d((\mu,w), (\mu',w'))= \delta(w(0), w'(0))+ D(\mu,\mu') \\
+ \int_0^{\langle
  \mu, 1 \rangle 
  \wedge \langle  \mu', 1 \rangle} dt\, \left(d_t(w_{\leq t}, w'_{\leq
    t})\wedge 1\right)+\bigl|\langle \mu,1\rangle-\langle \mu',1\rangle\bigr|, 
\end{multline}
where  $d_t$ is  the Skorohod  metric on  the space  $\D([0,t],  E)$ and
$w_{\leq t}$ denote the restriction  of $w$ to the interval
$[0,t]$. Notice $d$ is a distance on $\cw$. Indeed, we have that:
\begin{itemize}
   \item $d$ is
symmetric.
\item  $d((\mu,w), (\mu',w'))=0$ implies $D(\mu,\mu')=0$, that is
  $\mu=\mu'$, and then $w=w'$ a.e. on $[0,\langle \mu,1 \rangle)$. 
\item $d$ satisfies the triangular inequality. We  have for $(\mu,w),
  (\mu', w')$ and $(\mu'', w'')\in \cw$: 
\begin{align*}
  d((\mu,w), (\mu',w'))
&\leq  \delta(w(0), w''(0))+ \delta(w''(0), w'(0))+D(\mu,\mu'')
+D(\mu'',\mu')\\
&\hspace{1cm} +\int_0^{\langle
  \mu, 1 \rangle 
  \wedge \langle  \mu', 1 \rangle} dt\, \left(d_t(w_{\leq t}, w''_{\leq
    t})\wedge 1\right)\\
&\hspace{1cm} 
+\int_0^{\langle\mu, 1 \rangle 
  \wedge \langle  \mu', 1 \rangle} dt\, \left(d_t(w''_{\leq t}, w'_{\leq
    t})\wedge 1\right)+\bigl|\langle \mu,1\rangle-\langle
\mu',1\rangle\bigr|\\
&\leq  d((\mu,w), (\mu'',w''))+ d((\mu'',w''), (\mu',w'))\\
&\hspace{1cm} 
+ \left(\langle\mu, 1 \rangle   \wedge \langle  \mu', 1 \rangle 
-\langle\mu, 1 \rangle   \wedge \langle  \mu'', 1 \rangle \right)_+\\
&\hspace{1cm} 
+ \left(\langle\mu, 1 \rangle   \wedge \langle  \mu', 1 \rangle 
-\langle\mu'', 1 \rangle   \wedge \langle  \mu', 1 \rangle \right)_+\\
&\hspace{1cm} 
+\bigl|\langle \mu,1\rangle-\langle
\mu',1\rangle\bigr|\\
&\leq  d((\mu,w), (\mu'',w''))+ d((\mu'',w''), (\mu',w')),
\end{align*}
 \end{itemize} 
since
\[
(a\wedge b -a\wedge c)_++(a\wedge b -c\wedge b)_++|a-b|\leq |a-c|+|b-c|, 
\]
where $(x)_+=\max(x,0)$. 

We check  that $(\cw,  d) $ is  complete. Consider $((\mu_n,  w_n), n\in
\N)$  a Cauchy  sequence  in  $(\cw, d)$.  Since  $(\cm_f(\R_+), D)$  is
complete,  we  get that  $\mu_n$  converges to  a  limit  say $\mu$.  If
$\mu=0$,  the result  is clear.  If not,  for any  $\varepsilon>0$ small
enough, for $n$  and $n'$ large enough so  that $\langle \mu_n,1 \rangle
\wedge   \langle   \mu_{n'},1    \rangle\geq   \langle   \mu,1   \rangle
-\varepsilon$    we    deduce    from   \reff{eq:dist_d}    that,    for
$t_\varepsilon=\langle \mu,1  \rangle - 2\varepsilon$ and
$n_\varepsilon$ large enough,  $({w_n}_{\leq
  t_\varepsilon}, n\geq n_\varepsilon)$ 
is a Cauchy sequence in $\D([0,t_\varepsilon), E)$ and hence converge to
a   limit   $w_{\leq  t_\varepsilon}$.   Since   this   holds  for   any
$\varepsilon>0$  small enough,  we deduce  that $w$  is well  defined on
$[0,\langle  \mu,1 \rangle)$  and that  $w_n$  converges to  $w$ on  any
$\D([0,t),  E)$ for  $t<\langle  \mu,1 \rangle$.  We  deduce again  from
\reff{eq:dist_d} that $((\mu_n, w_n), n\in \N)$ converges to $(\mu, w)$.

We check that $(\cw, d)$ is separable. Let $(\mu_n, n\in \N)$ a dense
subset of  $(\cm_f(\R_+), D)$, and for each $n$, let $(w_{n,m}, m\in
\N)$ a dense subset of $(\D([0, \langle \mu_n, 1 \rangle), E),
d_{\langle \mu_n, 1 \rangle})$. Then, it is easy to check that $((\mu_n,
w_{n,m}); n,m\in \N)$ is a dense subset of $(\cw, d)$. 

Thus  the
space  $(\cw, d)$  is a  Polish space.   

We shall  write  $\mu_{\bar w}$
instead of $\mu$ when $\bar w=(\mu,w)$.
Recall \reff{def:H}. 
We consider a  family of probability measures $\bar \Pi_{x, \mu}$, for
$x\in E$ and the mass measure $\mu\in \cm_f(\R_+)$ on $\cw_x$, s.t. 
\begin{enumerate}[a)]
   \item $\mu_{\bar w}=\mu$, $\bar \Pi_{x,\mu}(d\bar w)$-a.s.;
   \item $w(0)=x$, $\bar \Pi_{x,\mu}(d\bar w)$-a.s.;
   \item $w$ has no fixed discontinuity: for all $s\in [0, \langle
     \mu,1 \rangle)$, $\bar \Pi_{x,\mu}(w(s-)=w(s))=1$;
   \item If $H(\mu)<\infty  $, then $w(\langle \mu,1  \rangle-)$ exists
     $\bar \Pi_{x,\mu}(d\bar w)$-a.s.; 
   \item If $H(\mu)<\infty  $ and $\nu\in \cm_f(\R_+)$, then
     under  $\bar \Pi_{x,  [\mu,\nu]}$,   $(w(r),  r\in  [0,  \langle  \mu,1
     \rangle)$  is  distributed  as  $(w(r),  r\in  [0,  \langle  \mu,1
     \rangle)$ under  $\bar \Pi_{x,\mu}$ and, conditionally  on $(w(r), r\in
     [0, \langle \mu,1 \rangle)$, $(w(r+\langle \mu,1 \rangle) , r\in
     [0,  \langle \nu,1 \rangle)$  is distributed  as $(w(r),  r\in [0,
     \langle \nu,1  \rangle))$ under $\bar \Pi_{w(\langle  \mu,1 \rangle-)
       ,\nu}$.
\end{enumerate}
The last property corresponds to the Markov property
conditionally on the mass measure. We shall assume that the mapping $(x,
\mu) \mapsto \bar \Pi_{x, \mu}$ is measurable. 

Let $\rho$  be an exploration process  starting at $\mu$.   We set $Y_t=
\langle \rho_t, 1 \rangle$. Recall  that $(Y_t, t\geq 0)$ is distributed
as a L\'evy process with Laplace exponent $\psi$ started at $\langle \mu,1
\rangle$. For $0\leq s<t$, we set $J_{s,t}=\inf_{s\leq u\leq t} Y_t$ and
$\rho_{s,t}=k_{(Y_s - J_{s,t})}  \rho_s=k_{(Y_t - J_{s,t})} \rho_t$, the
last equality being a consequence of the construction of the exploration
process. We also define $ \bar \rho_t^{(s)}$ as the unique measure $\nu$
s.t. $[\rho_{s,t}, \nu]=\rho_t$.

Conditionally on  $\rho$, we define a  probability transition
semi-group $R^\rho_{s,t}$
on  $\cw_x$  as follows:  for  $0\leq  s<t$ s.t. $J_{s,t}<Y_s$ or
$w(\langle \rho_s, 1 \rangle-)$ exists  and $\mu_{\bar w}=\rho_s$,  under
$R^\rho_{s,t} (\bar w,d\bar w')$ we have
\begin{enumerate}[i)]
   \item $\mu_{\bar w'}=\rho_t$, 
   \item $(w'(r), r\in [0, \langle \rho_{s,t}, 1
\rangle))=(w(r), r\in [0, \langle \rho_{s,t}, 1
\rangle))$, 
   \item $(w'(r), r\in [\langle \rho_{s,t}, 1
\rangle, \langle \rho_t, 1 \rangle))$ is distributed according to
$\bar \Pi_{w(  \langle \rho_{s,t}, 1 \rangle -) , \bar \rho_t^{(s)}}$. 
\end{enumerate} 
In (iii), by convention, if $\rho_{s,t}=0$, then $w( \langle \rho_{s,t},
1 \rangle -)  =x$.  Notice that for fixed  $s<t$, a.s.  $J_{s,t}<Y_s$ so
that, with the previous convention  $w( \langle \rho_{s,t}, 1 \rangle -)
$ is a.s. well defined.  Notice  that if $(\rho_s, w)$ is distributed as
$\bar  \Pi_{x, \rho_s}$,  then $(\rho_t,  w')$ is  distributed  as $\bar
\Pi_{x, \rho_t}$ thanks to condition e)  on $\bar \Pi$.  Thus we can use
the  Kolmogorov   extension  theorem  to  get that  there   exists  a  unique
probability   measure  $\P_{(\mu,w)}$   on  $(\cw_x)^{\R_+}$   s.t.  for
$0=s_0<s_1<\cdots <s_n$,
\begin{multline*}
   \P_{(\mu,w)}(W'_{s_0}\in A_0, \rho_{s_0}\in B_0,  \ldots,
W'_{s_n}\in A_n, \rho_{s_n}\in B_n)\\
=\E_\mu\left[\ind_{\{\rho_{s_0}\in B_0,  \ldots, \rho_{s_n}\in B_n\}}
  \ind_{\{w\in A_0\}}
\int_{A_1\times \cdots\times A_n}\!\!\!\!\!\!\! R^\rho_{s_0,s_1}(w,
dw_{s_1}) 
\cdots R^\rho_{s_{n-1},s_n}(w_{s_{n-1}}, dw_{s_n})\right]. 
\end{multline*}
We set $\bar W_s=(\rho_s, W'_s)$. 
Notice that $W'_s(r)=W'_t(r)$ for $r\in [0, \langle \rho_{s,t}, 1
\rangle)$ and thus that 
\[
d(\bar W_s, \bar W_t)\leq  D(\rho_s, \rho_t) + |Y_s \wedge Y_t -J_{s,t}|. 
\]
Since $\rho$ and  $Y$ are $\P_\mu$-a.s.  c\`ad-l\`ag, this  implies that the
mapping  $s\mapsto \bar W_s$  is $\P_{(\mu,w)}  $-a.s.  c\`ad-l\`ag  on $[0,
\infty )\bigcap \Q$.   Hence there is a unique  c\`ad-l\`ag extension to the
positive  real line,  we  shall  still denote  by  $\P_{(\mu,w)} $.  The
process   $(\bar   W_s,s\geq   0)$    is   under   $\P_{(\mu,w)}   $   a
time-homogeneous Markov process  living in $\D(\R_+, \cm_f(\R_+)\times
\cw)$.  
We call this distribution the distribution of the weighted L\'evy snake
associated with $\bar \Pi$. 

We  set $\cm_f^0(\R_+)$  the  set  of $\mu  \in  \cm_f(\R_+)$ such  that
$\supp(\mu)=[0,  H(\mu)]$  if   $H(\mu)<+\infty  $  and  $\supp(\mu)=[0,
H(\mu)) $ if  $H(\mu)=+\infty $. We the define $\Theta_x$  as the set of
all pairs $(\mu, w)\in \cw$  such that $\mu\in \cm_f^0(\R_+) $, $w(0)=x$
and at least one of the following three properties hold:
\begin{enumerate}
   \item[(i)] $\mu( H(\mu))=0$;
  \item[(ii)] $w(\langle \mu,1 \rangle-)$ exists;
\item[(iii)] $H(\mu)=+\infty $. 
\end{enumerate}

We denote  by $(\cf_s, s\geq  0)$ the canonical filtration  on $\D(\R_+,
\cm_f(\R_+)  \times   \cw)$.  One  can  readily  adapt   the  proofs  of
Propositions 4.1.1  and 4.1.2 of \cite{dlg:rtlpsbp} to  get the following
result.
\begin{theo}
   \label{theo:annexe}
   The process  $(\bar W_s, s\geq 0; \P_{(\mu,w)},  (\mu,w) \in \Theta_x
   )$ is a c\`ad-l\`ag Markov process in $\Theta_x$ and is strong Markov
   with respect to the filtration $(\cf_{s+}, s\geq 0)$.
\end{theo}

Let us remark that, when the family of probability measures
$\bar\Pi_{x,\mu}$ is just the law of a homogeneous Markov process
$\xi$ starting at $x$ and stopped at time $\langle\mu,1\rangle$, the
previous construction gives a snake with spatial motion $\xi$ and
lifetime process $X-I$, which is the total mass of the exploration
process. Notice that in \cite{dlg:rtlpsbp} the lifetime process is given
by the height of the exploration process. 

\subsection{The general L\'evy snake}

However, we need some dependency between the spatial motion and the
exploration process $\rho$ in order to  recover the usual  L\'evy snake
from  the weighted L\'evy snake.  Informally, we keep the spatial motion
from moving when time $t$ is ``on a mass'' of $\rho_s$. This idea can
be compared to a subordination and has already been used in the snake
framework by Bertoin, Le Gall and Le Jan in \cite{blglj:sbps} in order to
construct a kind of L\'evy snake from the usual Brownian snake.

Let $ \Pi_x$ be  the distribution of $\xi$ a c\`ad-l\`ag Markov
process taking values in $E$  with no fixed discontinuities and starting
at $x$,  such that  the mapping $x\mapsto  \Pi_x$ is  measurable. Recall
\reff{eq:def-ka}  and set $\hat  \mu_r=k_{\langle \mu,1  \rangle -r}\mu$
for $r\in [0,  \langle \mu,1 \rangle)$. We define  $\bar \Pi_{x,\mu}$ as
the distribution  of $(\mu,w)$ with  $w=(\xi_{H(\hat \mu_r)} ,  r\in [0,
\langle    \mu,1    \rangle))   $    under    $\Pi_x$.    Notice    that
$\xi_{r'}=w(\langle   \mu,  \ind_{[0,r']}   \rangle)$  for   $r'\in  [0,
H(\mu))$. In particular, $w$ is on constant on intervals
$\Big(\mu([0,r)), \mu([0, r])\Big] $ which corresponds to the atoms of
$\mu$.

We have that $\bar \Pi$ satisfies condition a)-e).

Let  $((\rho_s, W'_s), s\geq  0)$ be  the corresponding  weighted L\'evy
snake.  For $s\geq 0$,  $r\geq 0$,  we set  $W_s(r)=W'_s(\langle \rho_s,
\ind_{[0,r]} \rangle)$.  When $H$ is continuous,  the process $((\rho_s,
W_s),  s\geq  0)$   is  the  L\'evy  snake  defined   in  Section  4  of
\cite{dlg:rtlpsbp}  with  underlying   motion  $\xi$.   As a
consequence  of Theorem  \ref{theo:annexe},  we get  that the  (general)
L\'evy snake is strong Markov.

\begin{prop}
   \label{prop:annexe}
   The  process  $((\rho_s, W_s),  s\geq  0;  \P_{(\mu,w)}, (\mu,w)  \in
   \Theta_x  )$ is  a c\`ad-l\`ag  Markov process  in $\Theta_x$  and is
   strong Markov with respect to the filtration $(\cf_{s+}, s\geq 0)$.
\end{prop}

\subsection*{Acknowledgments}  The  authors  would  like to  thank  the
referee and the  associate editor for their comments, which in
particular helped to clarify the proof of the special Markov property.

\newcommand{\sortnoop}[1]{}

\end{document}